\documentclass[letterpaper, 11pt]{article}

\usepackage[margin=1in]{geometry}
\usepackage{setspace}\onehalfspace
\usepackage{microtype}

\usepackage{amsmath} \allowdisplaybreaks
\usepackage{amssymb}
\usepackage{amsthm}
\usepackage{array}

\usepackage{natbib}
\usepackage{bibentry}

\usepackage{booktabs} %
\usepackage{multirow}
\usepackage{multicol}
\usepackage{adjustbox}

\usepackage{graphicx}
\usepackage{appendix}
\usepackage[counterclockwise]{rotating} %
\usepackage{subcaption} %
\usepackage[normalem]{ulem} %
\usepackage{showexpl}
\usepackage{threeparttable} %
\usepackage{placeins}
\usepackage[hyphens]{url}
\usepackage[bookmarks=false,hidelinks]{hyperref}

\usepackage{tikz}
\usetikzlibrary{automata,calc,trees,positioning,arrows,chains,shapes.geometric,%
decorations.pathreplacing,decorations.pathmorphing,shapes,%
matrix,shapes.symbols,plotmarks,decorations.markings,shadows}

\usepackage{pgf}
\pgfdeclarelayer{background}
\pgfdeclarelayer{foreground}
\pgfsetlayers{background,main,foreground}

\usepackage[ruled,vlined, linesnumbered]{algorithm2e}
\usepackage{algpseudocode}
\algnewcommand\algorithmicforeach{\textbf{for each}}
\algdef{S}[FOR]{ForEach}[1]{\algorithmicforeach\ #1\ \algorithmicdo}

\usepackage{authblk}

\newcommand{\email}[1]{{\href{mailto:#1}{\nolinkurl{#1}}}}

\usepackage{bm}
\renewcommand{\vec}[1]{\boldsymbol{\mathbf{#1}}}
\newcommand{\mat}[1]{\vec{#1}}

\DeclareMathSymbol{\Gamma}{\mathord}{operators}{"00}
\DeclareMathSymbol{\Delta}{\mathord}{operators}{"01}
\DeclareMathSymbol{\Theta}{\mathord}{operators}{"02}
\DeclareMathSymbol{\Lambda}{\mathord}{operators}{"03}
\DeclareMathSymbol{\Xi}{\mathord}{operators}{"04}
\DeclareMathSymbol{\Pi}{\mathord}{operators}{"05}
\DeclareMathSymbol{\Sigma}{\mathord}{operators}{"06}
\DeclareMathSymbol{\Upsilon}{\mathord}{operators}{"07}
\DeclareMathSymbol{\Phi}{\mathord}{operators}{"08}
\DeclareMathSymbol{\Psi}{\mathord}{operators}{"09}
\DeclareMathSymbol{\Omega}{\mathord}{operators}{"0A}

\usepackage{xcolor}

\usepackage[pagewise]{lineno}
\usepackage{etoolbox}
\makeatletter
\def\do#1{\@namedef{#1c}{\ensuremath{\mathcal{#1}}}}
\docsvlist{A,B,C,D,E,F,G,H,I,J,K,L,M,N,O,P,Q,R,S,T,U,V,W,X,Y,Z}
\makeatother

\renewcommand{\bar}[1]{\mkern 1.5mu\overline{\mkern-1.5mu#1\mkern-1.5mu}\mkern 1.5mu}

\usepackage{mathtools}

\renewcommand{\vec}[1]{\boldsymbol{\mathbf{#1}}}

\def\Nc{\mathcal{N}}

\def\Gc{\mathcal{G}}

\def\Vc{\mathcal{V}}

\def\tr{\mathrm{tr}}
\def\dr{\mathrm{dr}}
\def\mask{\mathrm{mask}}

\renewcommand{\bar}[1]{\mkern 1.5mu\overline{\mkern-1.5mu#1\mkern-1.5mu}\mkern 1.5mu}
\renewcommand{\hat}{\widehat}
\renewcommand{\tilde}{\widetilde}

\usepackage{fancyvrb,listings}
\title{A Deep Reinforcement Learning Approach for Solving the Traveling Salesman Problem with Drone}
\author[1]{Aigerim Bogyrbayeva\thanks{\texttt{aigerim.bogyrbayeva@sdu.edu.kz}}}
\author[2]{Taehyun Yoon\thanks{\texttt{thyoon@unist.ac.kr}}}
\author[2]{Hanbum Ko\thanks{\texttt{hanbum.ko95@unist.ac.kr}}}
\author[3]{Sungbin Lim\thanks{\texttt{sungbin@korea.ac.kr}}}
\author[4]{Hyokun Yun\thanks{\texttt{yunhyoku@amazon.com}; This work does not relate to his position at Amazon.}}
\author[5]{Changhyun Kwon\thanks{\texttt{chkwon@usf.edu}; corresponding author. Part of this research was done while visiting KAIST, Daejeon, South Korea.}}
\affil[1]{Suleyman Demirel University, Kazakhstan}
\affil[2]{UNIST, South Korea}
\affil[3]{Korea University, South Korea}
\affil[4]{Amazon, U.S.A.}
\affil[5]{University of South Florida, U.S.A.}

\date{November 30, 2022}

\begin{document}
\maketitle

\begin{abstract}
Reinforcement learning has recently shown promise in learning quality solutions in many combinatorial optimization problems. In particular, the attention-based encoder-decoder models show high effectiveness on various routing problems, including the Traveling Salesman Problem (TSP). Unfortunately, they perform poorly for the TSP with Drone (TSP-D), requiring routing a heterogeneous fleet of vehicles in coordination---a truck and a drone. In TSP-D, the two vehicles are moving in tandem and may need to wait at a node for the other vehicle to join. State-less attention-based decoder fails to make such coordination between vehicles. We propose a hybrid model that uses an attention encoder and a Long Short-Term Memory (LSTM) network decoder, in which the decoder's hidden state can represent the sequence of actions made. We empirically demonstrate that such a hybrid model improves upon a purely attention-based model for both solution quality and computational efficiency. Our experiments on the min-max Capacitated Vehicle Routing Problem (mmCVRP) also confirm that the hybrid model is more suitable for the coordinated routing of multiple vehicles than the attention-based model. The proposed model demonstrates comparable results as the operations research baseline methods.
\end{abstract}

\section{Introduction}

Last-mile delivery, which refers to the transportation of products from a distribution center to the doorstep of a customer, is an integral part of the supply chain. 
However, last-mile delivery is often not cost-effective due to transportation costs associated with individualized shipments, complex routes, and various destinations. For instance, last-mile delivery costs comprise 50\% of the total delivery costs \citep{joerss2016parcel}. 
Therefore, emerging technologies such as unarmed aerial vehicles (drones) are viewed as a potential solution in reducing inefficiencies in last-mile delivery. 
For instance, Amazon launched the Prime Air program that aims to deliver parcels using drones \citep{amazon_news}. 
Similarly, Wing, a subsidiary of Alphabet, delivers books, meals, and medicine using drones \citep{wing_news}. 
The steady development of drone technology, including its capacity and flying range, enables the transport of various products. 
For instance, the HorseFly drone, used by UPS, can carry packages up to 10 pounds and has a flight time of 30 minutes \citep{UPS}.  

While drones certainly have high speed and do not require any infrastructure such as roads, bridges, etc., to fly, they also possess some limitations. Drones have a limited flying range, which depends on battery life. Also, drones cannot carry parcels of many customers; thus, they must return every time to a distribution center to pick up customer orders. On the other hand, trucks, which are traditionally used for last-mile delivery, have a large capacity to store all customer orders and can transport goods long distances. However, trucks usually have slow speeds due to congestion and require built infrastructure to reach their destinations. Therefore, combining the drone and truck is a promising tandem to improve the efficiency of last-mile delivery. Indeed, UPS has started testing combining the drone and the truck to deliver goods \citep{UPS}. In this setting, a truck driver loads a package into the drone and sends the drone to an autonomous route to an address. Meanwhile, the truck can serve other customers. When the drone returns to the truck, the driver swaps the battery of the drone and launches the next delivery. 

To effectively deploy trucks and drones together for last-mile delivery, we must answer several challenging questions such as which customers should be served by drone, which customers should be served by truck, where to recharge the drone, how to route both drone and truck, etc. In the literature, the problem of routing the truck and drone is known as the Traveling Salesman Problem with Drone (TSP-D) \citep{agatz2018optimization}. The objective of TSP-D is to serve customers in a minimal time, while the truck and the drone are subjected to routing constraints. TSP-D is an NP-hard problem \citep{agatz2018optimization}, highlighting the need to develop heuristic methods that can scale to large urban networks. The development of efficient computational methods for TSP-D 
\citep{agatz2018optimization,poikonen2019branch,roberti2021exact,vasquez2021exact} 
is still in its infancy. For instance, 
exact optimization approaches for TSP-D suffer from long computational time and are typically limited to smaller than 40-node problems.
Some heuristic optimization approaches can also take several minutes to a few hours to produce quality solutions for large-scale problems.
In this paper, we aim to develop a reinforcement learning approach for TSP-D, which is as competitive as fast heuristic optimization methods.

The main challenge of TSP-D stems from the \emph{simultaneous} decision-making for a heterogeneous fleet of vehicles and strong \emph{interdependency} among them. 
Unlike TSP or CVRP, where the decision is made for only one vehicle for a tour or multiple sequential tours, TSP-D requires route decisions for two different vehicles in coordination. 
The existing literature in reinforcement learning to route a single agent or a fleet of homogeneous vehicles in a multi-agent setting does not fit the nature of such heterogeneous vehicle routing problems. 
For instance, models to route a single truck or drone only handle cases when there is a single agent which interacts with an environment and does not generalize into the presence of several vehicles in the same environment. 
In contrast, in routing a heterogeneous fleet of vehicles, the interdependency among vehicles is a critical consideration since the capacity of one vehicle to fulfill customer orders is dependent on the interaction with the other vehicle. 
This, in effect, has a strong impact on the routing decisions of both vehicles.

Our goal is to devise an end-to-end learning model that can address the above-outlined challenges.
By end-to-end, we mean a learning-based model that takes problem instances as input and produces solutions as its direct output by itself without any additional algorithmic components.

Our main contributions are two-fold. 
First, we provide TSP-D as a Markov Decision Process (MDP) so that reinforcement learning can be developed and applied.
While MDP formulations for TSP only involve customer nodes, an MDP formulation for TSP-D should be able to capture the in-transit status of vehicles and remaining times to reach the next node, to express the state spaces of coordinated drone-truck routing adequately. Our MDP formulation is comprehensive and distinct from the dynamic programming (DP) approach of \citet{bouman2018dynamic}, which combines three DP problems sequentially to obtain the final formulation. This makes their DP formulation hard to use in end-to-end reinforcement learning approaches.

Second, we present a hybrid model (HM), based on an attention encoder and a Long Short-Term Memory (LSTM) decoder, for efficient routing of multiple vehicles taking into account required interactions between vehicles, along with a novel distributed training algorithm. 
Our HM consists of an attention-based encoder to encode a highly connected graph and an LSTM-based decoder that stores the routing history of all vehicles. 
Using a single decoder to route all vehicles allows passing information about the routing decisions of vehicles to each other, thus promoting efficient interaction. 
In training the hybrid model, we rely on a central controller, which observes an entire graph to route all vehicles. As our numerical results show, the proposed model performs comparably with operations research methods.

The remainder of the paper will proceed as follows: in Section \ref{sec:lit_rev} we present the literature review related to TSP-D and learning methods to solve routing problems; in Section \ref{sec:problem} we formally define TSP-D and present its Markov Decision Process formulation; in Section \ref{sec:methods} we present the Hybrid Model and a training algorithm to solve the problem; in Section \ref{sec:results} we present the main computational results for TSP-D; in Section \ref{sec:additional} we provide additional computational experiments in other routing problems; and finally in Section \ref{sec:concl} we conclude with final remarks and discussions.

\section{Literature Review} \label{sec:lit_rev}
In this section, we present recent studies aimed at solving TSP-D. 
We limit our literature review to studies that focus on routing a single drone and a single truck where both vehicles actively participate in serving customers. 
The detailed reviews of combined drone-truck routing and comparison among problem variants can be found in \citet{macrina2020drone} and \citet{chung2020optimization}.
Also, we present the recent advancements in learning methods to solve routing problems in general.

\subsection{Recent Studies to Solve TSP-D}

\citet{murray2015flying} first introduced the notion of \emph{sortie} for combined operations of the drone and truck, which they named the flying sidekick TSP (FSTSP). In a drone sortie, there are start and end nodes, which are defined as the nodes where the drone is launched and reunited with the truck. Between start and end nodes, the drone serves a single customer, while the truck may visit several customers. However, the start and end nodes of the drone cannot be the same, and they are limited to customer locations or a depot. The synchronization constraints are added to ensure that the drone and truck arrive at an end node at the same time. The drone routing is subjected to its battery life. At end nodes, the drone's battery is swapped with some service time for recharging. The objective of FSTSP is to minimize the total time spent in the system starting from when both the drone and truck leave a depot and return back after serving all customers. \citet{murray2015flying} presented the exact and heuristic methods to solve FSTSP. Later, \citet{ha2018min} extended the MIP of \citet{murray2015flying} to minimize the total transportation costs and \citet{yurek2018decomposition} introduced an iterative algorithm based on decomposition to solve FSTSP to minimize completion time. \cite{dell2021drone} and \cite{arxiv} present enhanced set of formulations of FSTSP, for which \cite{dell2021random} proposed a random restart local search heuristic to solve the large instances of the problem. \cite{dell2021algorithms} proposed a branch and bound algorithm to solve FSTPS up to 15 customers and devised an efficient heuristic to solve large instances. \cite{liu2022flying} used a reinforcement learning approach to consider stochastic traveling time in FSTSP. \cite{boccia2021column} develops an exact method relying on a column-and-row generation approach to solve FSTPS.

In contrast to drone sortie assumptions, \citet{agatz2018optimization} presented the idea of \emph{operations}, which is similar to drone sortie except it allows the drone to be launched and reunited with the truck at the same nodes. 
In fact, in their formulation, the truck and the drone are allowed to revisit nodes, which substantially increases the number of possible operations. 
Also, they considered different speeds for the drone and truck and allow the truck to wait for the drone at a rendezvous node.
With these problem settings, \citet{agatz2018optimization} define TSP-D.
\citet{carlsson2018coordinated} investigated the benefits of combined drone and truck operations in general, where drone launch and pick-up nodes are not limited to customer locations.
\citet{gonzalez2020truck} extended TSP-D by allowing drones to visit multiple customers per launch and present an iterative greedy search heuristic. 

For solving TSP-D, \citet{agatz2018optimization} devised a heuristic method that starts with an optimal TSP tour and partitions customer nodes into truck nodes and drone nodes using a dynamic programming approach. 
Later \citet{poikonen2019branch} introduced a branch-and-bound method, where the drone and truck can wait for each other at rendezvous nodes, and presented four heuristics to solve large instances of the problem. 
\citet{de2020variable} solved both TSP-D and FSTSP using a variable neighborhood search-based heuristic. 

Exact optimization methods are also developed for TSP-D.
\citet{bouman2018dynamic} developed a dynamic programming method that can solve instances with up to 15 customers to optimality, and \citet{poikonen2019branch} introduced a branch-and-bound method for up to 20 customer nodes. 
\cite{schermer2020b} developed a branch-and-cut approach to find optimal solutions for instances up to 20 customers. 
\cite{el2021parcel} presented a novel Mixed Integer Programming formulation for TSP-D, which was later improved using pre-processing and valid inequalities.
\citet{vasquez2021exact} proposed a Benders-type decomposition method and \citet{roberti2021exact} devised a branch-and-price algorithm to solve instances with up to 15 and 39 customer nodes, respectively.

\subsection{Learning Methods to Solve Routing Problems}
Recently reinforcement learning has received increased attention for solving routing problems.
Starting from \citet{bello2016neural}, which used Pointer Networks \citep{vinyals2015pointer} to solve TSP, learning methods have shown comparable results with optimization heuristics. 
For instance, \citet{dai2017learning} proposed a novel graph embedding structure to learn Q-function using DQN algorithm \citep{mnih2013playing} to solve a group of combinatorial optimization problems, including TSP. 

To directly learn the routing policies, the encoder-decoder structure has been a common choice for many studies. The encoder is a neural network with the primary goal to learn the graph representation given some input. The output of the encoder is passed to the decoder, a neural network that learns routing policy given the current state of the graph. Different neural network architectures have been proposed both for encoder and decoder, including fully attention-based models and LSTM-based models \citep{bogyrbayeva2022learning}. The attention mechanism is a technique that aims to find the most critical parts between inputs or between inputs and outputs. On the other hand, LSTM is a type of recurrent neural network that is designed to retain memories of both long and short-term dependencies in a sequence.

A series of research \citep{vinyals2015pointer, nazari2018reinforcement, kool2018attention} has shown that end-to-end learning, where only machine learning methods are applied, is a viable option for solving routing problems. 
In particular, \citet{kool2018attention} demonstrated that the Attention Model (AM), which leverages several layers of multi-head attention to learn representations of nodes, outperforms previous models, which rely on simpler representation models, for example, the model from \citet{nazari2018reinforcement} which used a single layer of single-head attention. 
This observation was consistent across several classes of routing problems and also with findings in natural language processing \citep{vaswani2017attention, devlin2018bert}. 

At each step of prediction, however, the decoder of AM conditions only on the current location of the vehicle and the current state of nodes and ignores the sequence of actions it has made in the past.
Such conditional independence is sensible for single-vehicle routing problems in which the ordering of past actions is irrelevant to future actions. 
However, it is not well-suited for problems like TSP-D requiring coordination among multiple vehicles.

Another set of studies focused on routing a fleet of vehicles consisting of homogeneous vehicles using a multi-agent model \citep{sykora2020multi,zhang2020multi} or a single-agent model \citep{bogyrbayeva2021reinforcement}. 
However, TSP-D considers routing multiple heterogeneous vehicles, including the truck and drone, with different routing constraints and capabilities. 

Beyond end-to-end learning models, there is another stream of research that combines optimization methods, multi-start approaches, or iterative processes with learning models to solve routing problems.
Examples include \citet{lu2019learning}, \citet{kwon2020pomo}, \citet{hottung2020neural}, \citet{kim2021learning}, and \citet{kool2021deep}.
However, all the above-mentioned papers do not focus on routing together heterogeneous vehicles such as drones and trucks, the defining feature of TSP-D.
Therefore, this paper focuses on developing an end-to-end model for TSP-D, which can provide a basis for integrated and iterative methods as well.

\section{Problem Statements and Formulations} \label{sec:problem}

In this section, we first define TSP-D and then present the Markov Decision Process formulation of the problem. 
\subsection{TSP-D Definition}

In TSP-D, we consider a complete graph $\Gc$ consisting of the set of nodes $\Nc = \{1,2,...,N\}$. Node 1 is a depot representing a warehouse, where a single truck equipped with a single drone is loaded with customer orders. 
Then, each of $N-1$ customers whose locations are known must be visited either by drone, truck, or both only once. 
When the distance from node $i$ to node $j$ is $d_{i,j}$, the traverse time from node $i$ to node $j$ for truck and drone is denoted by $\tau^\tr_{i, j}$ and $\tau^\dr_{i, j}$, respectively, where we assume $\tau^\dr_{i, j}\leq \tau^\tr_{i, j}$ for all $i, j \in \Nc$. 
Without loss of generality, we let the speed of truck and drone be 1 and $\alpha\geq 1$, respectively, so that $d_{i,j} = \tau^\tr_{i,j} = \alpha \tau^\dr_{i,j}$. We use the Euclidean distance to compute the traveling times of both the drone and the truck,
but our method straightforwardly generalizes to other distance metrics.

\begin{figure}
	\centering
	\tikzset{depot/.style={shape=circle, draw=black, fill=black, radius=0.35, text=white},
		trucknode/.style={shape=circle, draw=black, radius=0.35},
		dronenode/.style={shape=circle, draw=black, fill= gray, radius=0.35, text=white}
	}
	\begin{tikzpicture}[scale=0.95]
	\definecolor{green}{RGB}{31,182,83}
	\definecolor{red}{RGB}{203, 65, 84}
	\definecolor{blue}{RGB}{65, 105, 225}
	\definecolor{grey}{RGB}{192, 192, 192}
	\node[depot](depot) at (0, 0) {1};
	\node[left] at (depot.west) {\begin{tabular}{l}\scriptsize{$ \vec{s}_{0}=(\{1\}, 1, 1, 0, 0)$} \\  \scriptsize{$\vec{a}_0=(2, 2)$} \end{tabular}};
	\node[right] at (depot.east) {\begin{tabular}{l}\scriptsize{$\vec{s}_{6}=(\{1, 2, 3, 4, 5, 6, 1\}, 1, 1, 0, 0)$} \\ \scriptsize{$\vec{a}_6=(1, 1)$}\end{tabular}};
	\node[trucknode](truck1) at(0, 2) {2};
	\node[left] at (truck1.west) {\begin{tabular}{l}\scriptsize{$\vec{s}_{1}=(\{1, 2\}, 2, 2, 0, 0)$}\\ \scriptsize{$\vec{a}_1=(4, 3)$} \end{tabular}};
	\node[trucknode](truck2) at(0, 5) {4};
	\node[left] at (truck2.west) {\begin{tabular}{l}\scriptsize{$\vec{s}_{3}=(\{1, 2, 3, 4\}, 4, 4, 0, 0)$} \\\scriptsize{$\vec{a}_3=(6, 5)$} \end{tabular}};
	\node[trucknode](truck3) at(2.25, 3) {6};
	\node[right] at (truck3.east) {\begin{tabular}{l}\scriptsize{$\vec{s}_{5}=(\{1, 2, 3, 4, 5, 6\}, 6, 6, 0, 0)$} \\ \scriptsize{$\vec{a}_5=(1, 1)$} \end{tabular}};
	\node[dronenode](drone1) at(-2, 3.5) {3};
	\node[left] at (drone1.west) {\begin{tabular}{l}\scriptsize{$\vec{s}_{2}=(\{1, 2, 3\}, 4, 3, \tau^{\tr}_{2,4}-\tau^{\dr}_{2,3}, 0)$} \\  \scriptsize{$\vec{a}_2=(4, 4)$} \end{tabular}};
	\node[dronenode](drone2) at(2, 5) {5};
	\node[right] at (drone2.east) {\begin{tabular}{l}\scriptsize{$  \vec{s}_{4}=(\{1, 2, 3, 4, 5\}, 6, 5, \tau^{\tr}_{4,6}-\tau^{\dr}_{4,5}, 0)$} \\ \scriptsize{$\vec{a}_4=(6, 6)$} \end{tabular}};
	\draw[->] (depot) to (truck1);
	\draw[->] (truck1) to (truck2);
	\draw[->] (truck2) to (truck3);
	\draw[->] (truck3) to (depot);
	\draw[dashed, ->] (truck1) to (drone1);
	\draw[dashed, ->] (drone1) to (truck2);
	\draw[dashed, ->] (truck2) to (drone2);
	\draw[dashed, ->] (drone2) to (truck3);

	\end{tikzpicture}
	\caption{An example of TSP-D with state $\mathbf{s}_t := (\mathcal{V}_t, c_t^{\tr}, c_t^{\dr}, r_t^{\tr}, r_t^{\dr})$ and action $\mathbf{a}_{t}=(a_{t}^{\tr},a_{t}^{\dr})$ on six nodes: black, grey and white nodes represent a depot, customers served by drone and customers served by truck respectively. Solid and dashed lines correspond to drive and fly arcs, respectively. 
	The drone and the truck traverse together from the depot to the first customer. Then drone is launched to serve the second customer, while the truck traverses to serve the third customer. At this node, the drone reunites with the truck to be launched to serve the next customer. Drone and truck reunite at the last customer location and traverse back together to the depot.}
	\label{fig:tsp_example}
\end{figure}
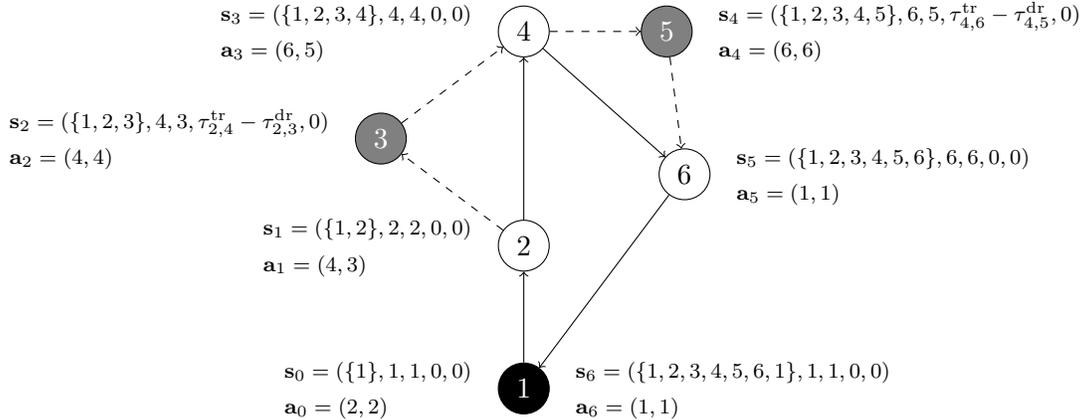

The objective of TSP-D is to complete serving all customers and return to the depot in a minimal time, while the truck and the drone are subjected to joint routing constraints.
For simplicity, we allow the drone to fly for an unlimited distance. 
However, the drone can serve only a single customer per launch, after which the drone must rendezvous with the truck to refresh its battery. 
Both the drone and truck are allowed to wait for each other at a rendezvous node. 
We assume that the customer service times, the parcel pick-up times, and the drone battery swapping times are negligible compared to the transportation time. 
Therefore, only customer nodes and the depot can be used to launch, recharge, and load the drone.
Figure~\ref{fig:tsp_example} shows a simple example and explains the MDP formulation proposed in the next section.

\subsection{MDP Formulation of TSPD}
To route both the drone and the truck in a shared urban environment, we assume that a central controller observes the entire graph and makes routing decisions for both the drone and the truck. We formalize this problem as a Markov Decision Process (MDP) with discrete steps and deterministic transitions. 

The state of the MDP is represented by the following variables. First, $\Vc_t \subseteq \Nc$ represents the set of nodes visited by any vehicle up to step $t$. 
We let $c_t^{\tr}$ and $c_t^{\dr}$ denote the current destination node of the truck and the drone, respectively. 
If the truck (or the drone) is in transit between two nodes at step $t$, then $c_t^{\tr}$ (or $c_t^{\dr}$) is set as the destination node of the vehicle. 
Or, if the truck (or the drone) is exactly located at a node, $c_t^{\tr}$ (or $c_t^{\dr}$) is set as the node index. 
By design, as we explain below, at least one vehicle will be located at a node at each step.
$r_t^\tr$ and $r_t^\dr$ represent the remaining times to arrive at destinations $c_t^{\tr}$ and $c_t^{\dr}$ respectively. 
When the truck (or the drone) is located at a node, $r_{t}^\tr = 0$ (or $r_{t} ^\dr = 0$). 
Then, $\vec{s}_t := (\Vc_t, c_t^{\tr}, c_t^{\dr}, r_t^{\tr}, r_t^{\dr})$ represents the state.
At step 0, vehicles are located at node 1. Therefore, $\Vc_0 = \{1\}$, $c_0^{\tr} = 1$, $c_0^{\dr} = 1$, $r_t^{\tr} = 0$, $r_t^{\dr} = 0$ and $\vec{s}_{0}=(\{1\}, 1, 1, 0, 0)$.

At each step, the central controller determines the next destination of the truck $a^{\tr}_t \in \Nc$ and that of the drone $a^{\dr}_t \in \Nc$. We cannot change the destination of a vehicle that is in transit from one node to another. When the drone is in transit, i.e. $r_t^\dr > 0$, we restrict $a^{\dr}_t = c_t^\dr$. Analogously, $a^{\tr}_t = c_t^\tr$ if the truck is in transit. For convenience, we refer to $\vec{a}_t := (a^{\tr}_t, a^{\dr}_t)$ as the action at step $t$.

Now let us discuss the state dynamics.
The next step occurs when either of the vehicles arrives at a node. 
Let $\hat{r}^\tr_t = r_t^\tr + \tau_{c_t^{\tr}, a_t^{\tr}}$
and
$\hat{r}^\dr_t = r_t^\dr + \tau_{c_t^{\dr}, a_t^{\dr}}$ be 
``updated'' versions of $r^\tr_t$ and $r^\dr_t$ right after the action is selected,
with the convention $\tau_{i,i} = 0$ for any $i \in \Nc$.
Then, the time taken until the next step, which happens when at least one vehicle arrives at the destination, will be the minimum of the two. We denote this as $C_t := \min(\hat{r}^\tr_t, \hat{r}^\dr_t)$. 
The remaining time is updated as $r_{t+1}^\tr = \hat{r}_t^\tr - C_t$, 
$r_{t+1}^\dr = \hat{r}_t^\dr - C_t$. The set of visited nodes is updated as $\Vc_{t+1} = \Vc_t \cup \{ a_t^{v} \mid  r_{t+1}^v = 0, v \in \{\tr, \dr \} \}$. Also, $c_{t+1}^{v} = a_t^v$ for $v \in \{\tr, \dr \}$.

Let $T$ be the index of the step when both the drone and the truck return back to the depot after serving all customers. Technically, we can define the value of $T$ as infinity if the controller fails to serve all customers and return vehicles. Then, our cost function (negative reward) is the total time spent in the system or the makespan:
\begin{align}
    C = \sum_{t=0}^{T} C_t \label{eq:cost}.
\end{align}

\section{Solution Methods} \label{sec:methods}

In this section, we discuss in detail the proposed Attention-LSTM Hybrid Model and present training methods. 

\subsection{The Attention-LSTM Hybrid Model}
Following \citet{vinyals2015pointer}, we aim to learn a probabilistic routing policy $\pi_{\theta}$ as a product of conditional probabilities
\begin{align}
	\pi_{\theta}(\vec{a}_1, \vec{a}_2, \ldots,\vec{a}_{T} |\Gc) = \prod_{t=1}^{T}p_{\theta}(\vec{a}_{t}|\vec{s}_{t}),
	\label{eq:policy}
\end{align}
where $p_{\theta}$ is a function parameterized by $\theta$. Then, we employ policy gradient methods \citep{williams1992simple} to learn $\theta$. In practice, encoder-decoder architecture has shown its effectiveness in solving routing problems \citep{bogyrbayeva2022learning}. Therefore, we also utilize an encoder-decoder structure that enables learning both a graph representation through encoding and policy $p_{\theta}$ through decoding.

\begin{figure}[t]
	\centering 
	\includegraphics[width=\columnwidth]{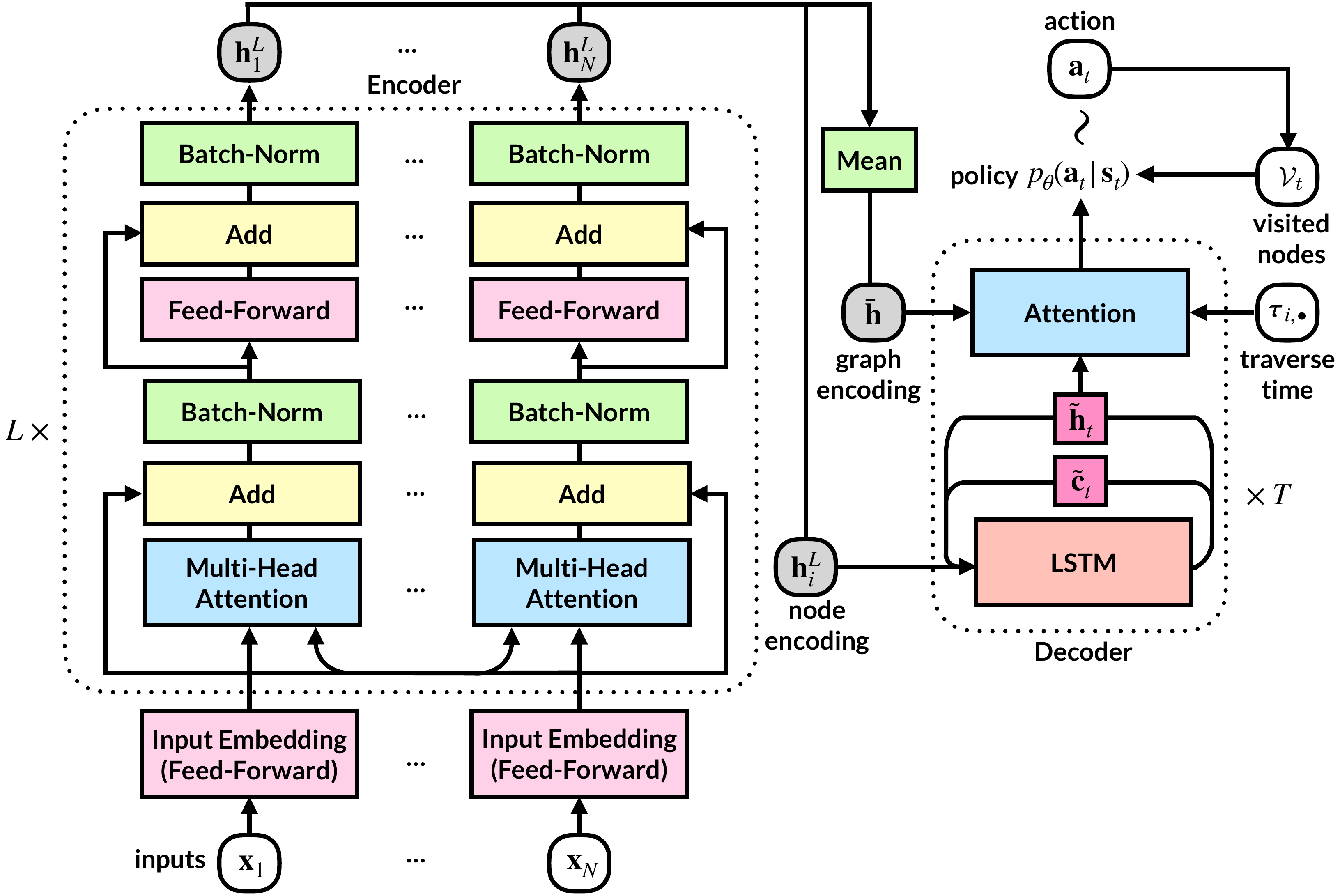}
	\caption{An Encoder-Decoder structure of Attention-LSTM Hybrid model (HM).}
	\label{fig:hybrid}
\end{figure}

In particular, to coordinate multiple vehicles with different capacities, we propose a hybrid model (HM) that uses multi-head attention for encoding as in \citet{kool2018attention}, but leverages LSTM hidden states for decoding in order to address dependencies between subsequent actions. 
Since, in TSP-D, the relative locations of the drone and the truck are a key piece of information for coordinated routing, LSTM's ability to store past decisions in its hidden states will be more effective than stateless attention-based decoders.
The model by \citet{nazari2018reinforcement} also utilizes an LSTM-based decoder for single-vehicle routing problems such as CVRP, but AM outperforms the \citeauthor{nazari2018reinforcement} Model (NM). 
We argue that our HM is suitable for coordinated routing by combining AM's encoder and NM's decoder with additional features. 
The overall structure of HM is shown in Figure~\ref{fig:hybrid}. 
We formally describe HM below.

\paragraph{Encoder}
Given a graph with a set of nodes $\Nc$, we have coordinates of each node in $\mathbb{R}^{2}$. Let $\vec{x}_n = (x_n, y_n)\in\mathbb{R}^{2}$ represent coordinates for node $n \in \Nc$. To embed such a fully connected graph, we start from the initial embeddings of the nodes  using linear transformation:
\begin{equation}
	\vec{h}^0_n = \mat{W}^0\vec{x}_n + \vec{b}^0 \quad \forall n \in \Nc,
\end{equation}
where $\mat{W}^0$ and $\vec{b}^0$ are learnable parameters. 
These initial embeddings of nodes $\{\vec{h}_{n}^{0}:n\in\mathcal{N}\}$ are then passed through $L$ number of attention layers. 
Each attention layer, $l$, consists of two sublayers: a multi-head attention layer and a fully connected feed-forward layer.

\emph{A multi-head attention (MHA) layer} takes as an input the output of the previous layer (either output from the initial embedding or the output of the previous attention layer) and passes messages between nodes. In particular, for each input to MHA we compute the values of $\vec{q} \in \mathbb{R}^{d_q}, \vec{k} \in \mathbb{R}^{d_q}, \vec{v} \in \mathbb{R}^{d_v} $ or queries, keys and values respectively by projecting the input $\vec{h_n}$: 
\begin{equation}
	\vec{q}_n = \mat{W}^Q\vec{h}_n, \quad 
	\vec{k}_n = \mat{W}^K\vec{k}_n, \quad 
	\vec{v}_n = \mat{W}^V\vec{h}_n \quad \forall n \in \Nc,
\end{equation}
where $\mat{W}^Q, \mat{W}^K, \mat{W}^V$ are trainable parameters with sizes ($d_k \times d_h$), ($d_k \times d_h$), and ($d_v \times d_h$), respectively. From keys and queries, we compute compatibility between nodes $i$ and $j$ as follows since we consider the fully connected graph:
\begin{equation}
	u_{i,j} = \dfrac{\vec{q}^\top_{i}\vec{k}_j}{\sqrt{d_k}} \quad \forall i,j \in \Nc .
\end{equation}
We use compatibility $u_{i,j}$ to compute the attention weights, $a_{i,j} \in [0, 1]$ using softmax:
\begin{equation}
	a_{i,j} = \dfrac{e^{u_{ij}}}{\sum_{j^{'}}{e^{u_{i,j^{'}}}}} \label{msg}.
\end{equation}
Then a message received by node $n$ is a convex combination of messages received from all nodes:
\begin{equation}
	\vec{h}^{'}_{n} = \sum_{j}a_{n,j}\vec{v}_j.
\end{equation}
Instead of using a single head, we use a multi-head attention layer with head size $M$, which allows passing different messages from other nodes. Then after computing messages from each head using \eqref{msg}, we can combine all massages coming to node $n$ as follows:
\begin{equation}
	\vec{MHA}_{n}(\vec{h}_1, \dots, \vec{h}_N) = \sum^M_{m=1}{\mat{W}^O_m\vec{h}^{'}_{nm}}.
\end{equation}
The output of the MHA sublayer, along with the skip connection, is passed through Batch Normalization:
\begin{equation}
	\vec{\hat{h}}^{l}_{n} = \vec{BN}^l(\vec{h}_{n}^{l-1}+\vec{MHA}_{n}(\vec{h}^{l-1}_1, \dots, \vec{h}^{l-1}_N)).
\end{equation}
The output of Batch Normalization is then passed through a fully connected feed-forward (FF) network with the ReLU activation function. We again apply the skip connection and batch normalization to the output of FF. 
\begin{equation}
	\vec{h}^l_{n} = \vec{BN}^l(\vec{\hat{h}}^{l}_n + \vec{FF}^l(\vec{\hat{h}}^{l}_n)),
\end{equation}
where 
\begin{equation}
	\vec{FF}^l(\vec{\hat{h}}^{l}_n)) = \mat{W}^{\text{ff}, 1}\cdot\text{ReLU}(W^{\text{ff}, 0}\vec{\hat{h}}^{l}_n + \vec{b}^{\text{ff},0}) + \vec{b}^{\text{ff},1}.
\end{equation}

\paragraph{Decoder}
Given the encoder outputs as embeddings of each node in the graph, we denote by $\bar{\vec{h}}$ the graph encoding computed as the mean of the embeddings of all the nodes in a graph. Then to select an action for a decision taker (either the drone or truck), we pass to the LSTM the embedding of the last node selected by the decision taker denoted as $\vec{h}^{L}_{i}$ where $i$ represents the current location: 
\begin{equation}
	\tilde{\vec{h}}'_{t+1}, \tilde{\vec{c}}'_{t+1} = \text{LSTM}(\vec{h}^{L}_{i}, (\tilde{\vec{h}}_{t}, \tilde{\vec{c}}_{t})).
\end{equation}
Here $\tilde{\vec{h}}_{t}$ and $\tilde{\vec{c}}_{t}$ correspond to the hidden and cell states at time $t$.
We apply dropout to the output of LSTM with probability $p$:
\begin{align}
    \tilde{\vec{h}}_{t+1} = \text{Dropout}(\tilde{\vec{h}}'_{t+1}, p),\quad \tilde{\vec{c}}_{t+1} = \text{Dropout}(\tilde{\vec{c}}'_{t+1}, p).
\end{align}
Then we pass the hidden state of the LSTM to attention along with the graph embedding and the traveling time from the current location of a decision taker, similar to \citet{bogyrbayeva2021reinforcement}, to speed up training and improve the performance. 
The attention vector is calculated as follows:
\begin{equation} \label{eq:att_vec}
	\vec{a}_{i, \boldsymbol{\cdot}} = \vec{v}_a^\top \tanh \bigg(\mat{W}^a [\bar{\vec{h}};  \vec{\tilde{h}}_{t+1}; \mat{W}^d \vec{\tau}_{i,\boldsymbol{\cdot}}] \bigg),
\end{equation}
where $i$ represents the current location of a decision taker, 
$\vec{\tau}_{i, \boldsymbol{\cdot}}$ is the vector of the traverse times from the current node to all other nodes in a graph, and $\vec{v}_a, \mat{W}^a,$ and $\mat{W}^d$ are trainable parameters with dimensions $d_h$. The resulting attention vector is $\vec{a}_{i, \boldsymbol{\cdot}}$ and $a_{i,j}$ denotes its element.
Then the attention is passed through softmax to produce the probabilities of visiting the next node as follows:
\begin{equation}
	p_{\theta}(a_{t}^{v}=j|\mathbf{s}_{t}) = \begin{cases}\dfrac{\exp(a_{i,j})}{\sum_{j^{'} \in \Nc \setminus \Vc_t}{\exp(a_{i,j^{'}})}} &: j \in \Nc \setminus \Vc_t, \\
	0 &: j \in \Vc_t,
	\end{cases}
\end{equation}
so that we do not allow nodes to be revisited by either vehicle. 
In general, optimal solutions to TSP-D can have certain nodes visited more than one time \citep{agatz2018optimization}. 
However, we empirically found from preliminary experiments that preventing vehicles from revisiting nodes improves the efficiency of training and does not degrade the quality of final solutions. 
Hence, we disallow revisiting in all our experiments. 
See Section~\ref{sec:revisiting} for further discussion on revisiting.

\subsection{Training Methods}

The presented HM must be trained through efficiently exploring actions and receiving feedback in a form of rewards. 

The proposed model determines the probability distribution, $\pi_\theta(\Vc_T|\Gc)$, which, given a graph, produces the sequence of nodes to be visited by the drone and truck. 
In TSP-D, we aim to minimize the makespan $C$ in \eqref{eq:cost}. 
We define the training objective function as follows:
\begin{align}
    J(\theta|\Gc) = \mathbb{E}_{\pi_{\theta}(\Vc_T|\Gc)}[(C - b(\Gc)) \log \pi_\theta(\Vc_T|\Gc)],
\end{align}
where $b(\Gc)$ is a baseline for variance reduction. 
We use a critic network to estimate $b(\mathcal{G})$, while an actor network is used to learn the policy $\pi_{\theta}$; hence we can compute the gradients $\nabla_{\theta}J(\theta|\mathcal{G})$ using the REINFORCE algorithm \citep{williams1992simple}. 
 
\paragraph{Distributed RL Training}
\begin{figure*}
\centering
	\includegraphics[width=0.8\columnwidth,keepaspectratio]{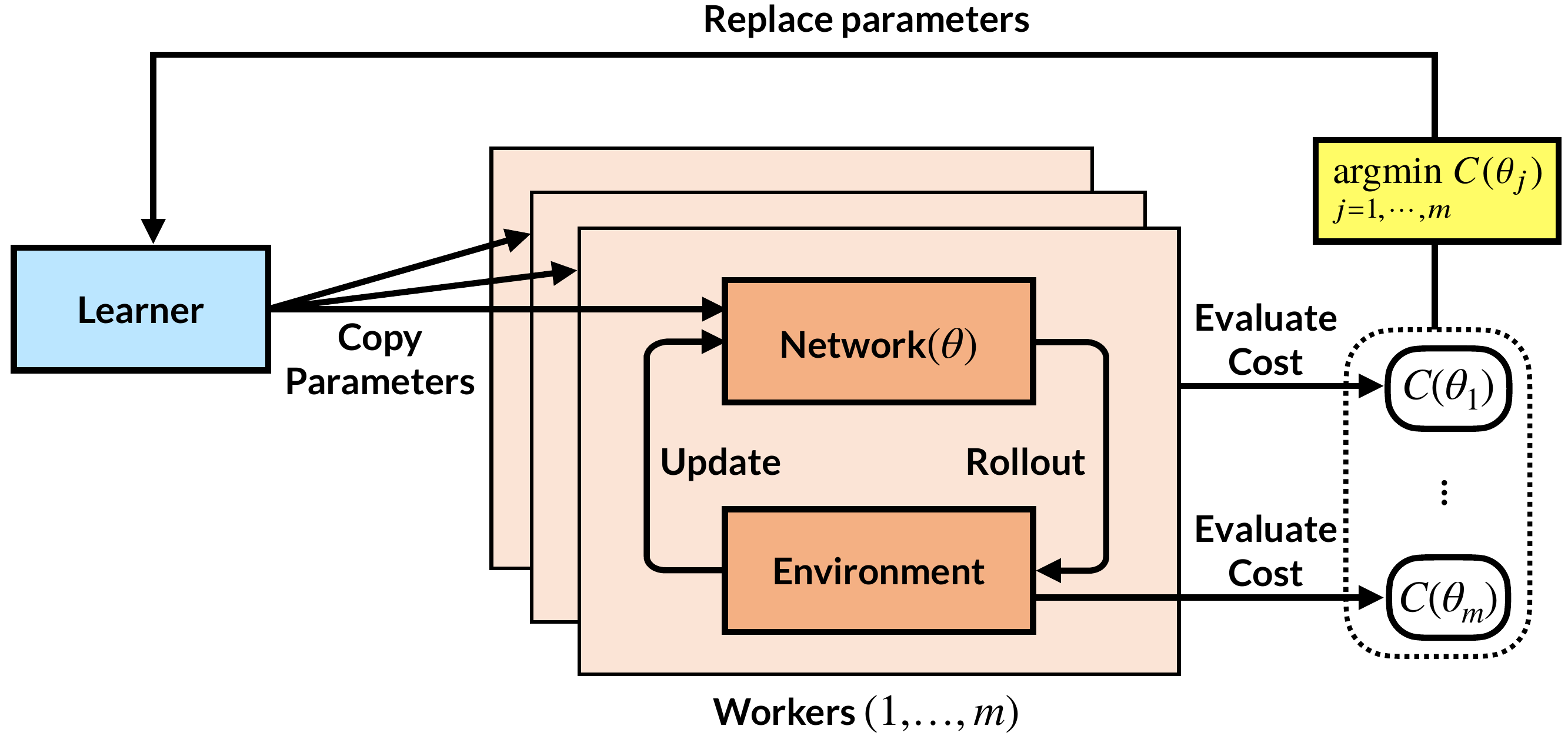}
	\caption{  
	Overview of distributed RL training for HM: each of multiple workers executes rollout with different training instances, updates its parameters, and evaluates costs. The learner replaces its network parameters with the ones from the best worker.}
	\label{fig:DFPG_pic}
\end{figure*}

To accelerate and stabilize the training of HM, we devise a new training algorithm that takes advantage of parallelization. 
We modify the Advantage Actor-Critic (A2C) algorithm to develop the \textit{Distribute-and-Follow Policy Gradient} (DFPG) training algorithm. 
To achieve scalable efficiency and reduced variance, DFPG selects gradients of the best performing worker when updating the central learner's weights as shown in Algorithm \ref{alg:TSP-D}. 
We generate multiple parallel workers of an RL agent with the same network parameters copied from the central learner. 
Each worker performs a rollout with a batch of episodes of different data instances. 
Then the network parameters of the workers are updated by gradient using the advantage function as a baseline (see Algorithm \ref{alg:TSP-D_rollout}). 
Next, we evaluate the updated network parameters of each worker with the shared validation instances. 
Finally, we select the worker whose performance is the best among the other workers in a greedy manner where the central learner copies the network parameters from the best worker as shown in Figure \ref{fig:DFPG_pic}.  

\begin{figure}[h!]
    \centering
	\includegraphics[width=0.45\columnwidth, keepaspectratio]{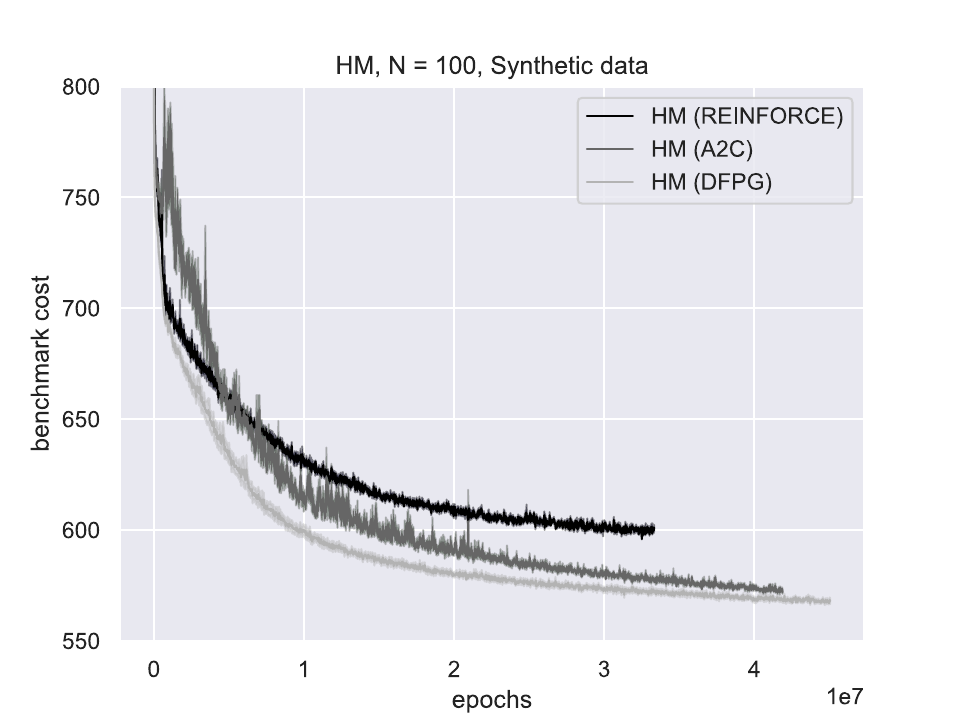}
	\includegraphics[width=0.45\columnwidth, keepaspectratio]{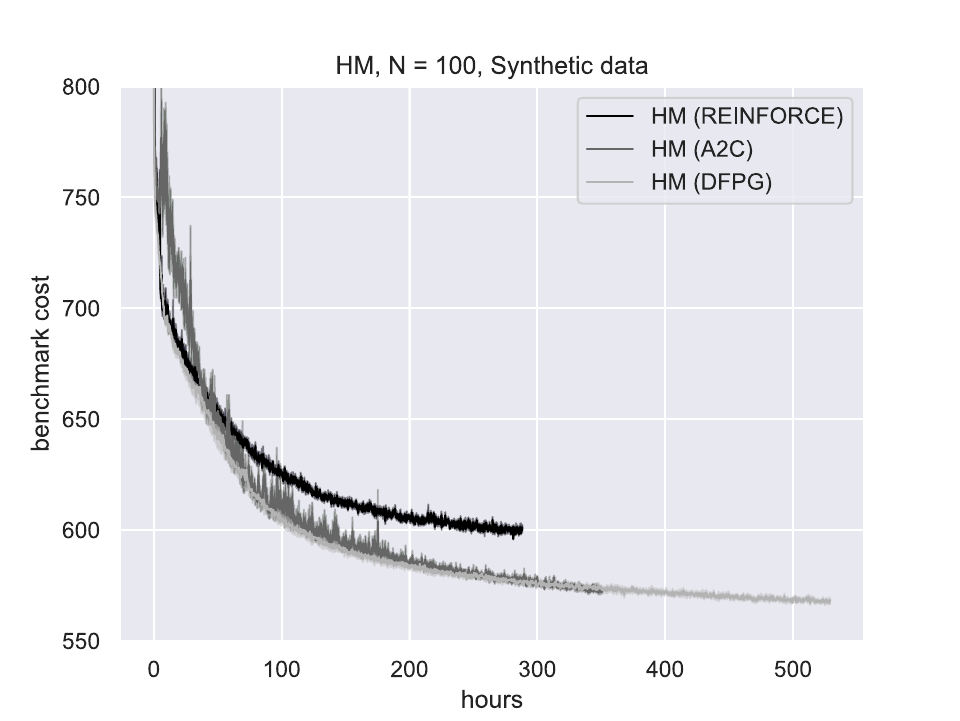}
	\caption{Benchmark cost curve comparison between A2C and DFPG on 100-node graphs from the environment without revisiting. Synthetic data refers to 100 benchmark instances used in Table \ref{tab:TSPD_uniform}.}
	\label{fig:env_no_revisit_n100}
\end{figure}

Our experiments show that DFPG stabilizes the final performance of the proposed HM at the convergence point, shown in Figure \ref{fig:env_no_revisit_n100}. 
We also observe that DFPG accelerates the training of AM, but only a slight improvement over AM trained with REINFORCE is shown at the final epoch. 
Hence, we apply the REINFORCE algorithm for training AM as suggested in \citet{kool2018attention} for a fair comparison. 

\begin{algorithm}[!hbpt]
	\caption{Distribute-and-Follow Policy Gradient}
	\label{alg:TSP-D}
	Generate a set of $K$ multiple parallel workers, $\Omega = \{\omega_1, \dots, \omega_K\}$.\\
	Initialize actor and critic networks parameters $\{\theta^a_1, \dots, \theta^a_K\}$, $\{\theta^c_1, \dots, \theta^c_K\}$ of $A$. \\
	Set the maximum number of epochs, $M_\mathrm{epochs}$ \\
	\For{$\mathrm{epochs} = 1$ \KwTo $M_\mathrm{epochs}$} 
	{
		Initialize an empty list of validation rewards $\textbf{R}$ \\
		\For(\Comment{in parallel}){$k$ = 1 \KwTo $K$} 
		{
			Initialize network gradients $d \theta^a_k \gets 0, d \theta^c_k \gets 0$ \\
			$B_k$ $\sim$ DataGenerator($\rho$) \\
			Call Rollout($\omega_k$, $B_k$) \\
			$R_k \gets \text{Evaluate}(\omega_k)$ \\
			Append $R_k$ to $\textbf{R}$ \\
	    }
	    $j \gets \arg \min\limits_{1 \leq j \leq K}(\textbf{R})$\\
	    \For(\Comment{in parallel}){$k$ = 1 \KwTo $K$}
	    {
	        $\theta^a_k \gets \theta^a_j$ \\
	        $\theta^c_k \gets \theta^c_j$ \\
	    }
    }
\end{algorithm}
\FloatBarrier

\begin{algorithm}[!htbp]
	\caption{Rollout}
	\label{alg:TSP-D_rollout}
	\KwIn{Batch of data $B$ with  a number of episodes denoted $M_{epis}$, a worker id $k$. Set the maximum number of steps denoted, $T$}
	Initialize reward $R$ \\
	Initialize LSTM initial state $(\tilde{h}_0, \tilde{c}_0)$\\
	$x_0$, $\mask_0 \gets \textsc{Env.Reset}$(B) \\
	\For{t=0 \KwTo $T$}{
	$a^\text{tr}_t, (\tilde{h}_{t'}, \tilde{c}_{t'}) \gets \pi_{\theta_k^a}(x^{\text{tr}}_t$, $\mask^{\text{tr}}_t$, $(\tilde{h}_t, \tilde{c}_t))$ \\
	$a^{\text{dr}}_t, (\tilde{h}_t, \tilde{c}_t) \gets \pi_{\theta_k^a}(x^{\text{dr}}_t$, $\mask^{\text{dr}}_t$, $(\tilde{h}_{t'}, \tilde{c}_{t'})$) \\
	$x_{t+1}$, $\mask_{t+1}$, $C_t \gets \textsc{Env.Step}(a^{\text{tr}}_t,a^{\text{dr}}_t $) \\
	$R \gets R - C_t$
	}
	Calculate $b(x_0; \theta^c_k)$ using critic
	
	$d \theta^a_k \gets \frac{1}{M_\text{epis}}\sum_{m=1}^{M_\text{epis}} (R^m-b^m(x^m_0; \theta^c_k))\nabla_{\theta^a_k}\log \pi_{\theta_k^a}$ \\
	$d \theta^c_k \gets \frac{1}{M_\text{epis}}\sum_{m=1}^{M_\text{epis}} \nabla_{\theta^c_k}(R^m-b^m(x^m_0; \theta^c_k))^2$
	
\end{algorithm}

\section{Computational Studies} \label{sec:results}
We empirically demonstrate the strength of HM against learning and optimization heuristics for TSP-D.
In all experiments, we assumed the drone is twice as fast as the truck; i.e. $\alpha=2$.

\subsection{Training and Evaluation Configurations}

\paragraph{Data Generation}
We generate two sets of datasets for TSP-D from different types of locations. In \emph{random locations} dataset, we randomly sample x and y coordinates of each node from a uniform distribution over $[1, 100] \times [1, 100]$, with the exception of the depot node which is distributed over $[0, 1] \times [0, 1]$. This makes the depot always located at the corner. 
This generation method is identical to \citet{agatz2018optimization}.

While such a uniform sampling is standard in the TSP-D literature, we create a new dataset named \emph{Amsterdam} for more realistic experiments, from the dataset originally used in \citet{haider2019optimizing}. In particular,
to identify potential demand locations for urban mobility and logistics services, we used an electric vehicle (EV) parking location dataset used in \citet{haider2019optimizing}. 
As EVs are usually parked on the urban streets, next to the curbside chargers, these locations reflect where potential customers are located. 
For a detailed description of the dataset, refer to \citet{haider2019optimizing}.
We randomly pick depot and customer nodes from the entire Amsterdam dataset to create various problem instances.
We use the Amsterdam dataset for evaluation only and generate the training data on the fly using kernel density estimation.
In particular, we estimate probability density functions of x and y coordinates of the customer nodes using Gaussian kernels. 
We deploy the Gaussian-KDE library of SciPy \citep{2020SciPy-NMeth} that selects bandwidth using Scott's Rule \citep{scott2015multivariate}. 

\paragraph{A critic network for HM} 
The architecture of the critic network has similarities to the actor network, except we do not use recurrent neural networks. 
The goal of the critic network is to estimate the total time needed to serve all the customers and return back to the depot, which is achieved through embedding the initial state of the graph. 
In particular, we embed x and y coordinates of the nodes through element-wise projections with 1D convolution networks whose outputs are passed through attention followed by feed-forward networks. 
We use Python 3.8 and PyTorch 1.5 for all our implementations.

\paragraph{Hyperparameters} We use 128 instances of data generated on the fly per epoch to train HM. 
We use 3 layers in the encoder and 8 attention heads. 
We initialize all the encoder parameters using Uniform(1/$\sqrt{d}$, 1/$\sqrt{d}$), with $d$ the input dimension. 
In the decoder, we initialize the LSTM hidden and cell states with zeros and apply dropout with $p=0.1$. 
Embedding dimensions are set to 256 for TSP-D with uniformly distributed graphs and 128 for other experiments. 
We use $d_h=128$ to compute attention in the decoder.
We use the constant learning rate set to $10^{-4}$. 
For AM, we used the same hyperparameters as in \citet{kool2018attention}.

\paragraph{Gap} In all the experiments, we use the following formula to report the relative gap between two solution methods:
\begin{align}\label{gap}
   \text{Gap}_i = \frac{Z_i- Z_i^*}{Z_i^*} \times 100 \%,
\end{align}
where $Z_i$ is the cost of a solution method for each instance $i$ and $Z_i^*$ is the cost of the best-performing solution method among all methods compared for instance $i$. Then we report `Gap' as the mean value of all $\text{Gap}_i$ values for the same setting.

\paragraph{Baselines}
Even though there is a surge of studies proposing the combination of learning methods with optimization heuristics to solve routing problems outlined in \citet{mazyavkina2021reinforcement}, to analyze our neural architecture choice for HM, we compare against AM \citep{kool2018attention}, which is a high-performing \emph{end-to-end} learning algorithm for routing problems. 
We also compare with NM \citep{nazari2018reinforcement}.
Since our HM is built upon AM's encoder and NM's decoder, they are natural choices for comparison.
We also compare against heuristic optimization algorithms from the operation research literature, namely `TSP-ep-all' \citep{agatz2018optimization} and `divide-and-conquer heuristic' (DCH) \citep{poikonen2019branch}. 
The TSP-ep-all heuristic starts with an initial TSP tour, then \emph{partitions} the TSP tour into nodes to be served by drone and truck, followed by local \emph{search} algorithms for the TSP tour and subsequent partitioning.
On the other hand, the DCH of \citet{poikonen2019branch} splits the initial TSP tour into several subgroups and partition nodes within each subgroup using a branch-and-bound algorithm. 
We create a new DCH approach by applying TSP-ep-all to partition each subgroup instead of branch-and-bound; we call this method the Divide-Partition-and-Search (DPS) algorithm.
In particular, we use DPS/$g$, which divides all nodes into subgroups so that each subgroup has $g$ nodes. 
While $g=10$ was originally used in \citet{poikonen2019branch}, we also test with $g=25$. 
As $g$ increases, the solution time and quality also increase. 
Note that when $g=N$, the DPS/$N$ is identical to TSP-ep-all. 
We implemented TSP-ep-all and DPS/$g$ in Julia 1.5 \citep{bezanson2017julia}, while the initial TSP tours are obtained by the Concorde TSP Solver \citep{applegate2001tsp}.

\paragraph{Decoding Strategies} There are two methods by which we sample solutions from the trained hybrid model. 
In the first method called \emph{greedy}, we always select nodes with the highest probabilities to visit at each time step. 
In the second method, called \emph{sampling}, we sample multiple solutions independently from the trained model according to \eqref{eq:policy}. 
Then, we select the minimal-cost sample as the solution.

\paragraph{Solution times} 
We measure the inference time of the benchmark dataset using an NVIDIA A100 GPU (80 GiB) and AMD EPYC 7713 64-Core Processor CPU (128 threads used) for HM and AM. 
For heuristic methods, Intel Xeon E5-2630 2.2 GHz CPU is used.
Each TSP-D benchmark dataset for each size consists of 100 instances.
For HM greedy, NM greedy, and AM greedy, instances are evaluated not in a batch but one by one in a sequence on the GPU. 
Heuristic methods use the same strategy but use a single thread of the CPU instead of a GPU.
For HM sampling, denoted by HM ($s$), we make $s$ copies of each instance and sample solutions for all $s$ copies in a batch on the GPU. 
AM sampling is executed in the same way as HM sampling. 
The entire code, including the environment for HM and AM, supports GPU parallelization, boosting the solution time up to 4.5 times from the environment with CPU computations for HM (4800) and AM (4800) on $N=100$ instances.
Note that we do not intend to compare the solutions times of the heuristics and the learning methods directly; rather, we aim to understand the trends.
We report the average time for a single instance.

\subsection{Results on TSP-D with Unlimited Flying Range}

\begin{figure}
    \centering
    \begin{tabular}{| >{\centering\arraybackslash} m{1.3cm} | >{\centering\arraybackslash} m{5.5cm} | >{\centering\arraybackslash} m{5.5cm} |}
    \hline
    Instance & Optimal Solution & HM Greedy  \\
    \hline
    1        
    & \includegraphics[width=5.3cm]{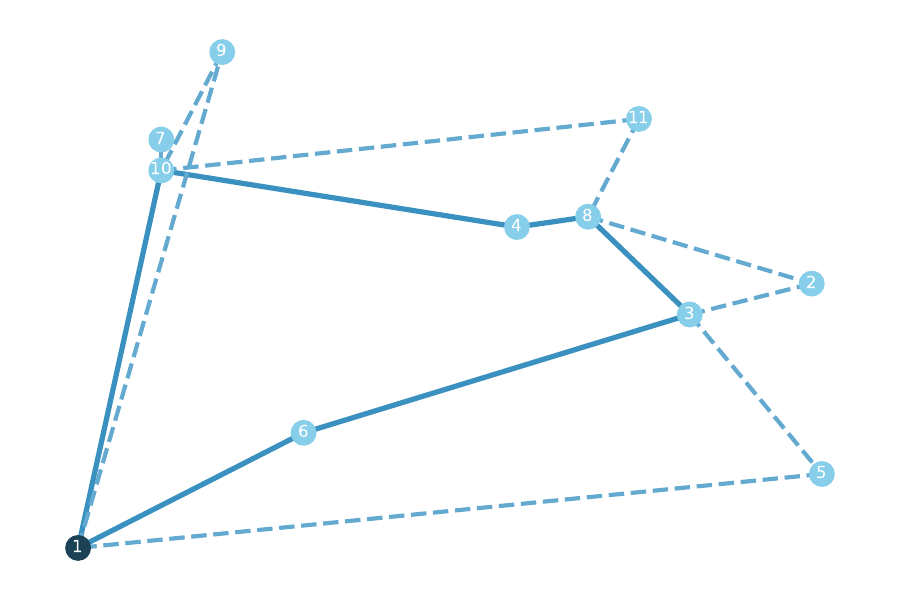}
    Cost = 221.19
    & \includegraphics[width=5.3cm]{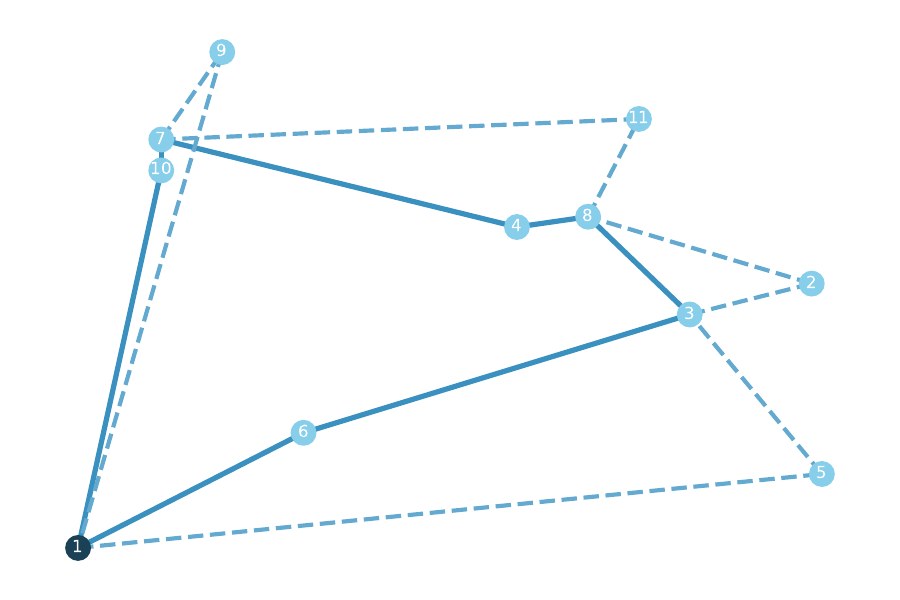}
    Cost = 223.41 (Gap = 1.00\%) \\
    \hline
    2        
    & \includegraphics[width=5.3cm]{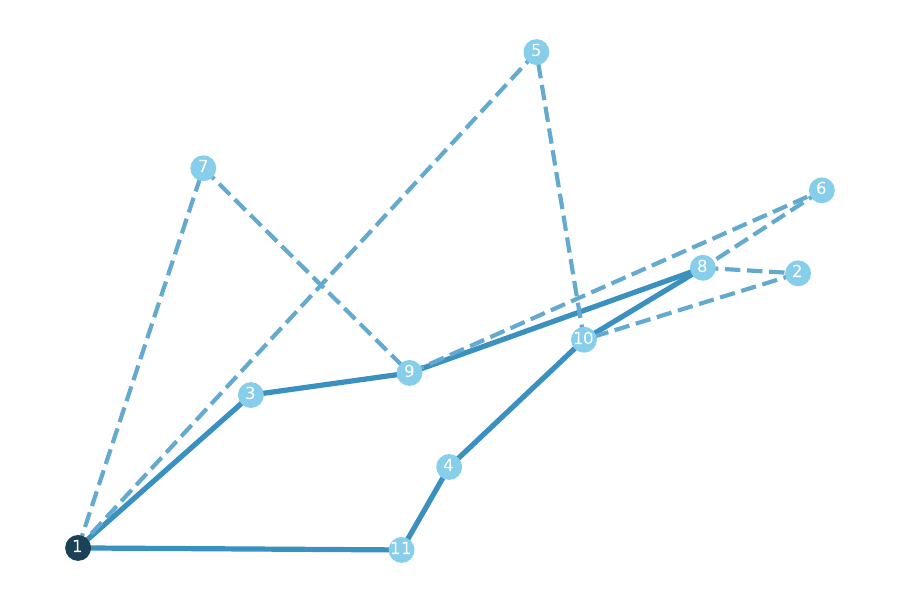}
    Cost = 205.76
    & \includegraphics[width=5.3cm]{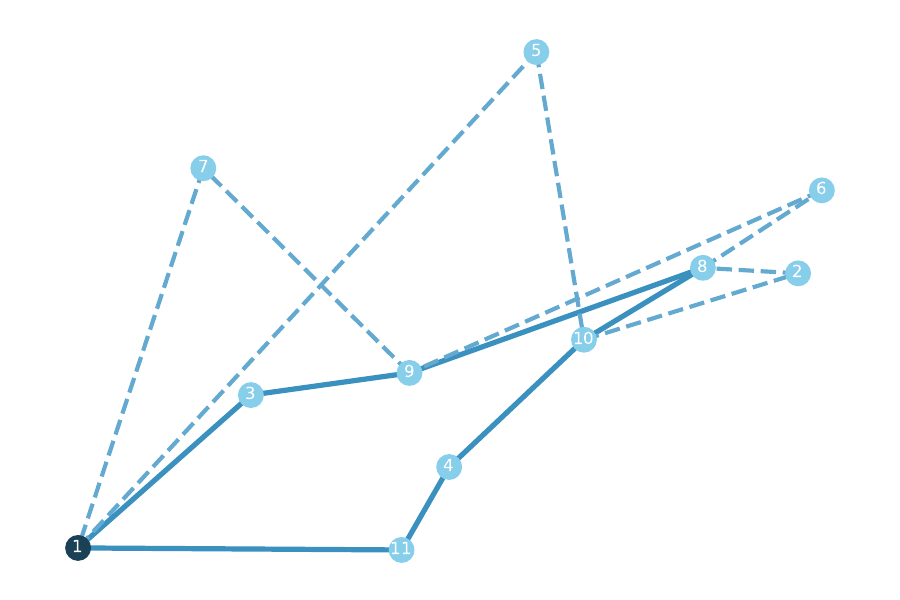}
    Cost = 205.76 (Gap = 0.00\%) \\
    \hline
    3       
    & \includegraphics[width=5.3cm]{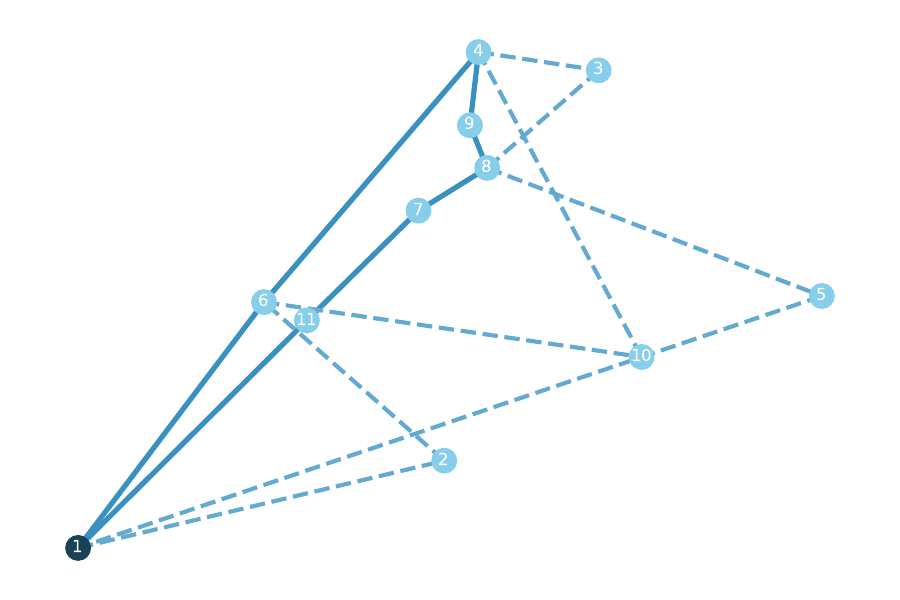}
    Cost = 192.96
    & \includegraphics[width=5.3cm]{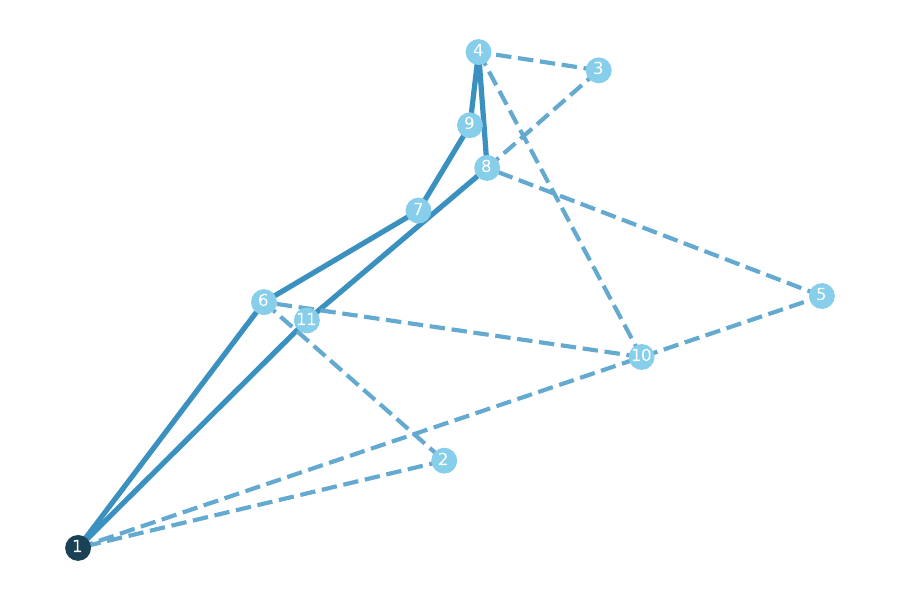}
    Cost = 193.99 (Gap = 0.53\%) \\
    \hline
    4       
    & \includegraphics[width=5.3cm]{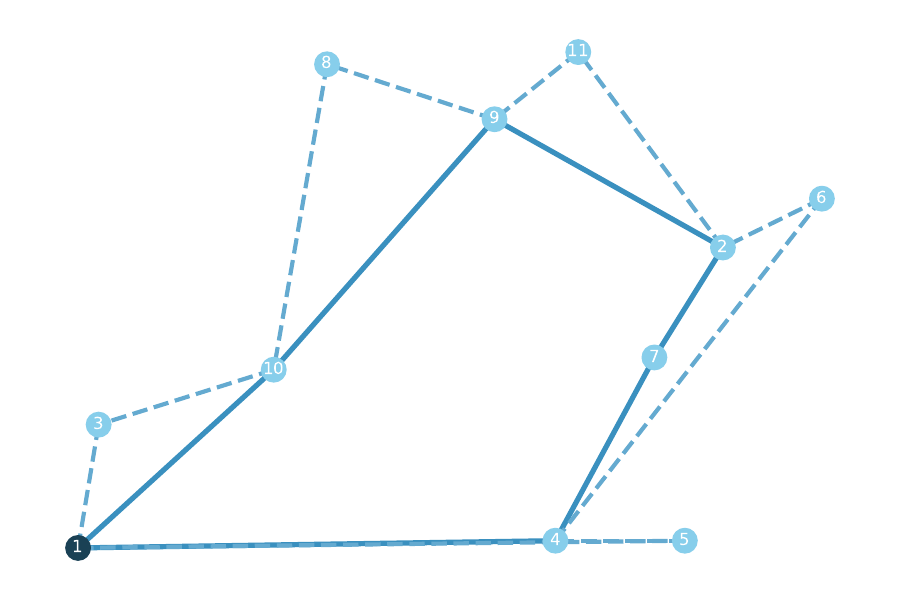}
    Cost = 241.26
    & \includegraphics[width=5.3cm]{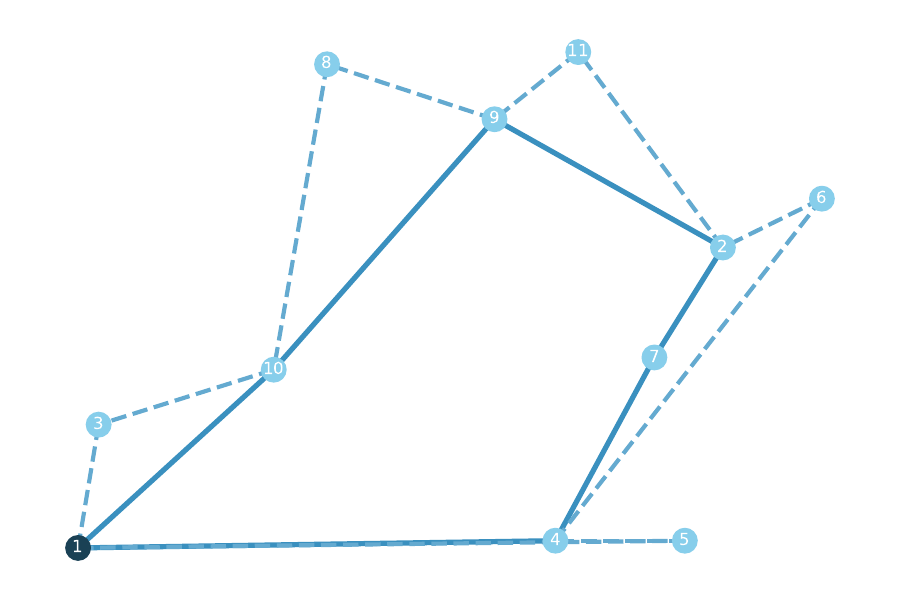}
    Cost = 241.26 (Gap = 0.00\%) \\
    \hline
    5        
    & \includegraphics[width=5.3cm]{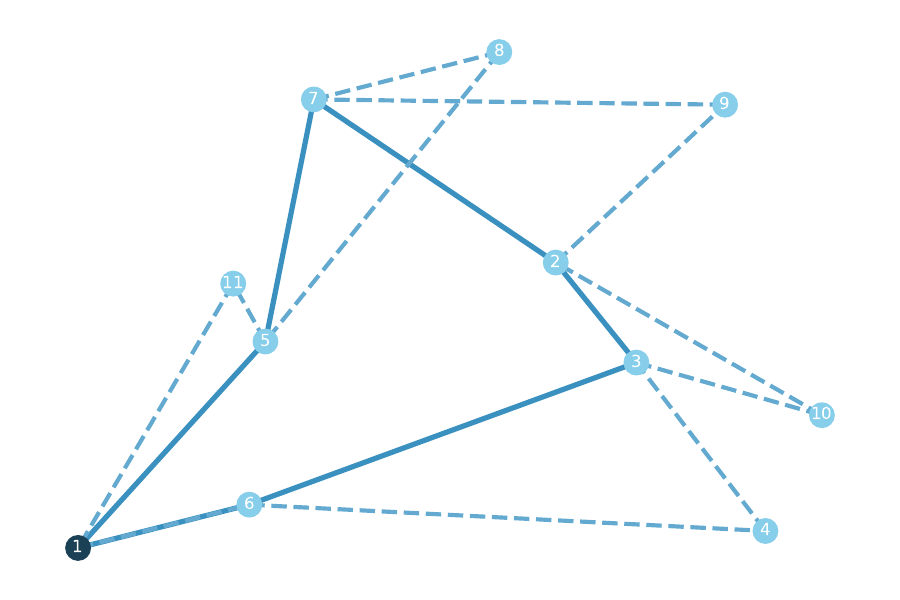}
    Cost = 248.14
    & \includegraphics[width=5.3cm]{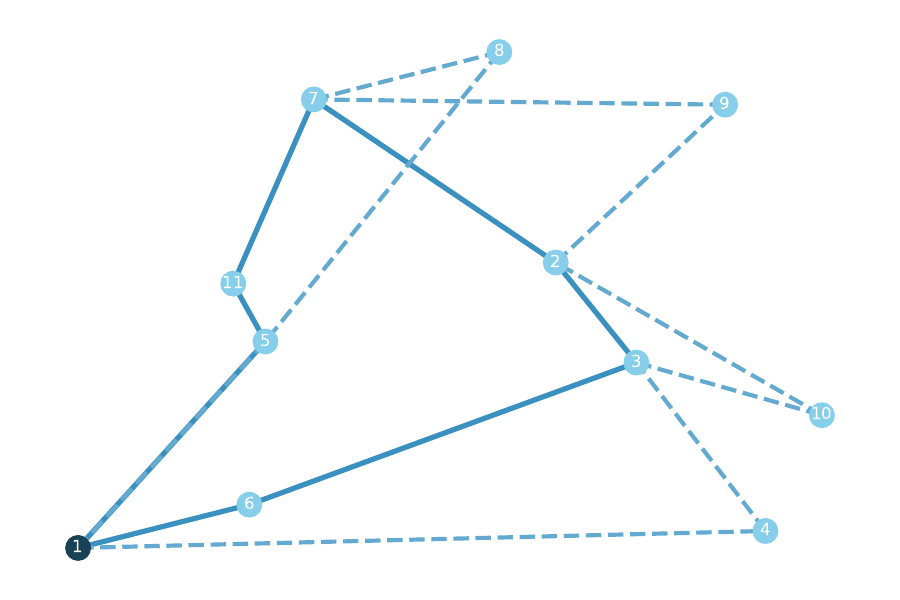}
    Cost = 249.85 (Gap = 0.69\%) \\
    \hline    
    \end{tabular}
    \caption{Routing examples: optimal vs. HM greedy (Instances 1 to 5). Node 1 is the depot, the solid line is the truck route, and the dashed line is the drone route.}%
    \label{fig:example1}
\end{figure}

\begin{figure}
    \centering
    \begin{tabular}{| >{\centering\arraybackslash} m{1.3cm} | >{\centering\arraybackslash} m{5.5cm} | >{\centering\arraybackslash} m{5.5cm} |}
    \hline
    Instance & Optimal Solution & HM Greedy  \\
    \hline
    6       
    & \includegraphics[width=5.3cm]{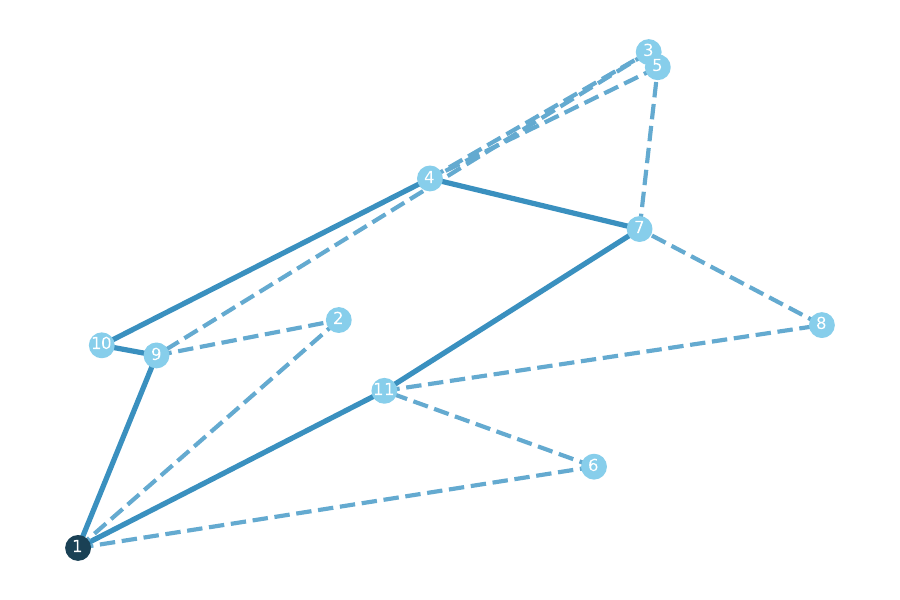}
    Cost = 217.69
    & \includegraphics[width=5.3cm]{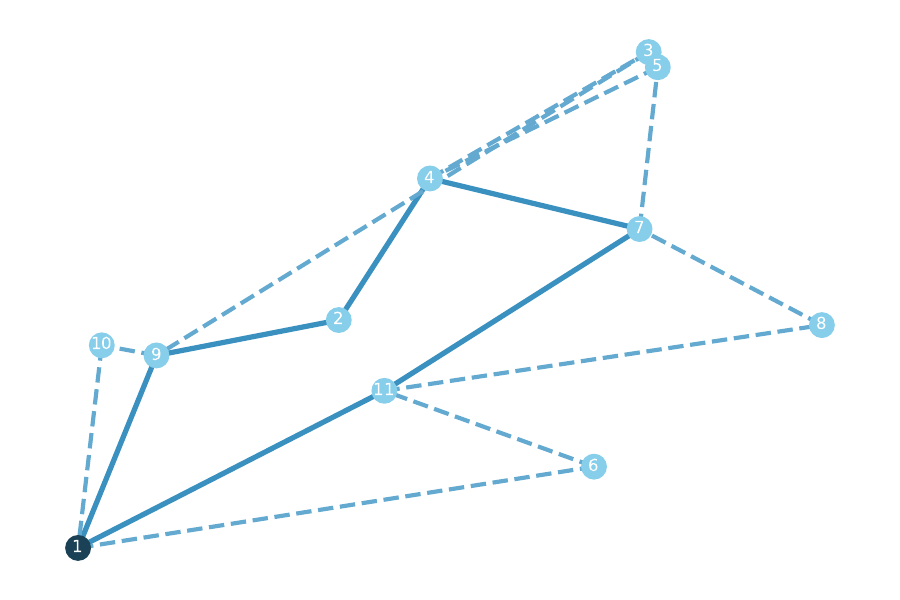}
    Cost = 217.69 (Gap = 0.00\%) \\
    \hline
    7        
    & \includegraphics[width=5.3cm]{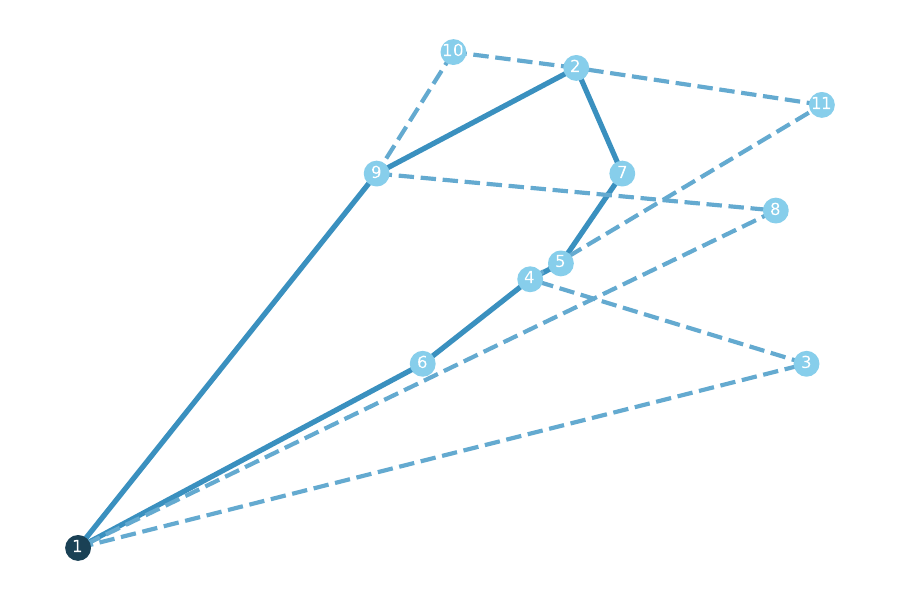}
    Cost = 237.34
    & \includegraphics[width=5.3cm]{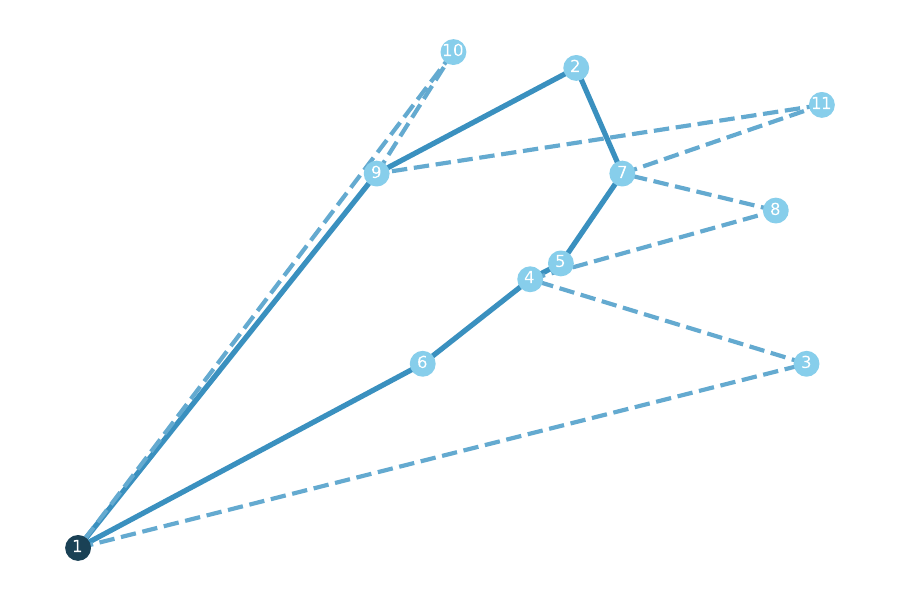}
    Cost = 240.47 (Gap = 1.32\%) \\
    \hline
    8       
    & \includegraphics[width=5.3cm]{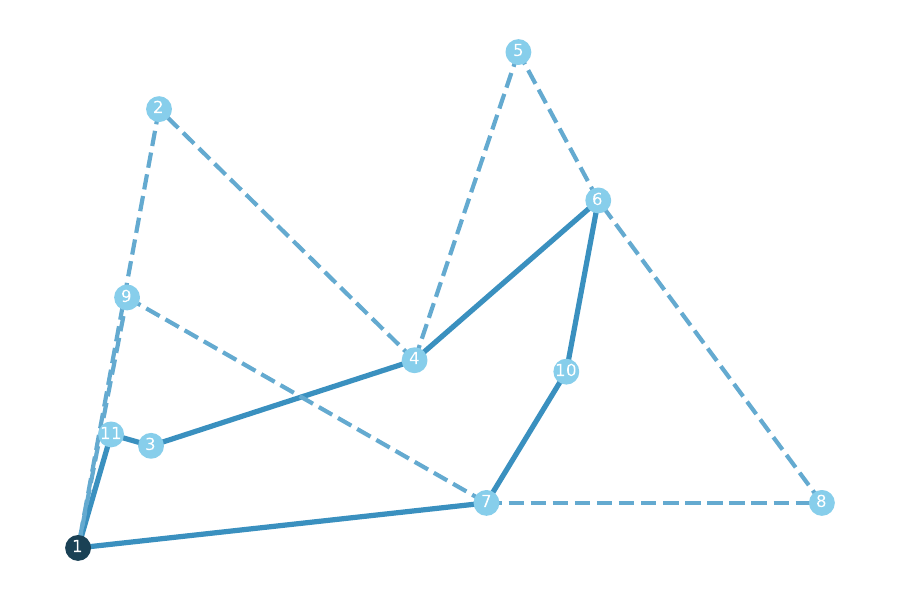}
    Cost = 214.77
    & \includegraphics[width=5.3cm]{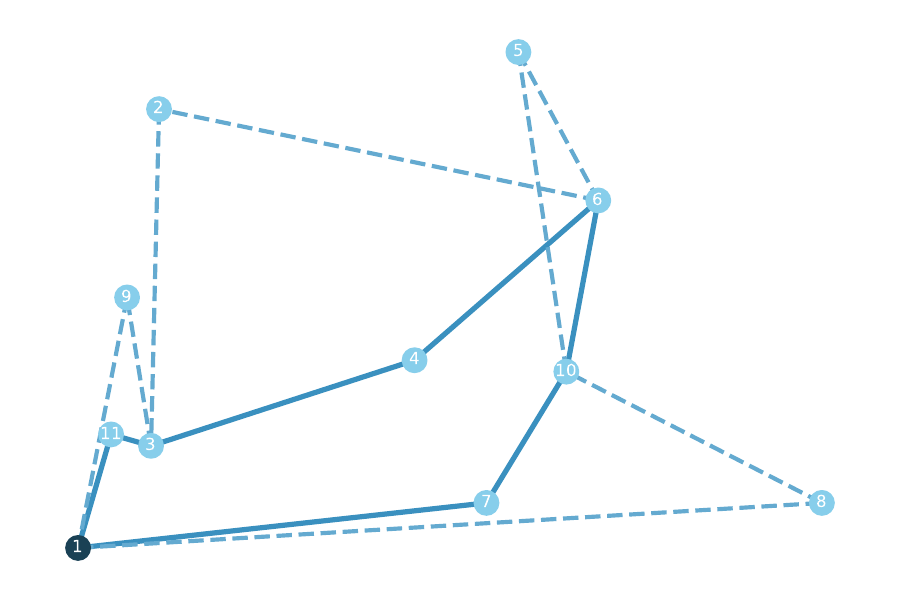}
    Cost = 226.64 (Gap = 5.53\%) \\
    \hline
    9       
    & \includegraphics[width=5.3cm]{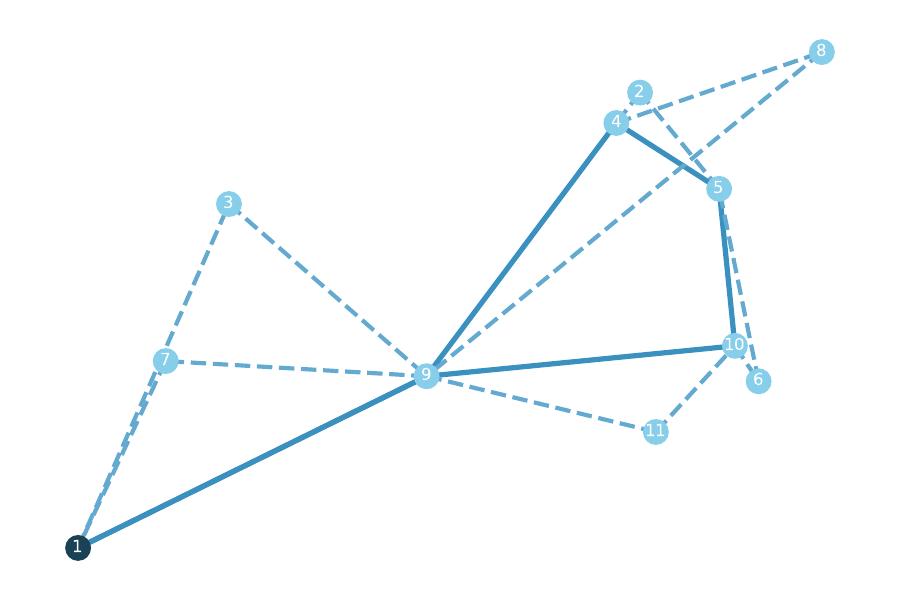}
    Cost = 256.34
    & \includegraphics[width=5.3cm]{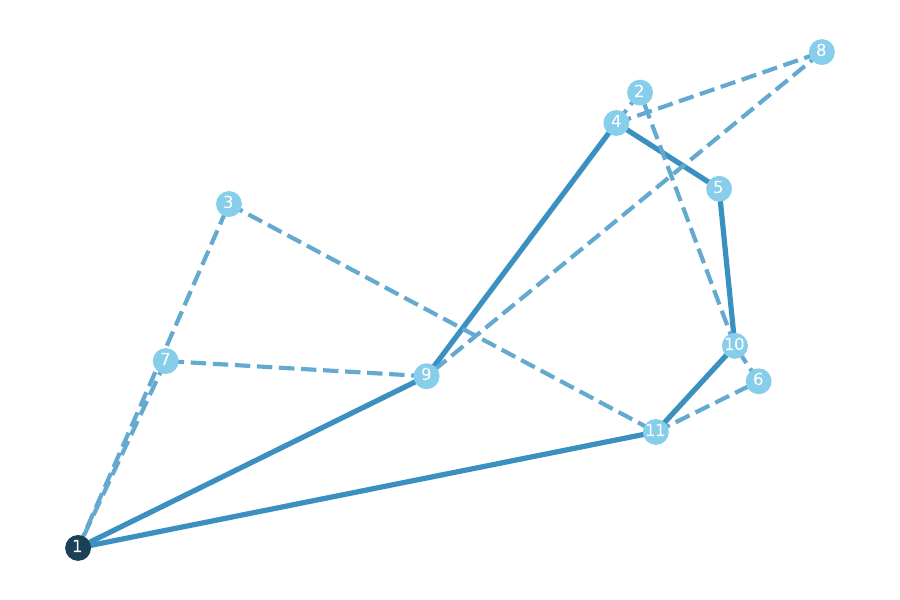}
    Cost = 256.83 (Gap = 0.19\%) \\
    \hline
    10        
    & \includegraphics[width=5.3cm]{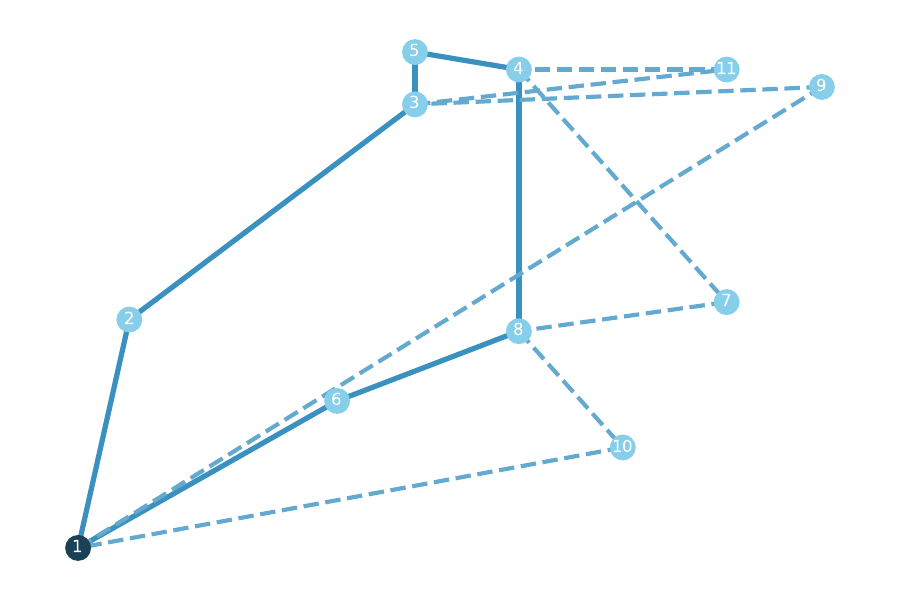}
    Cost = 227.90
    & \includegraphics[width=5.3cm]{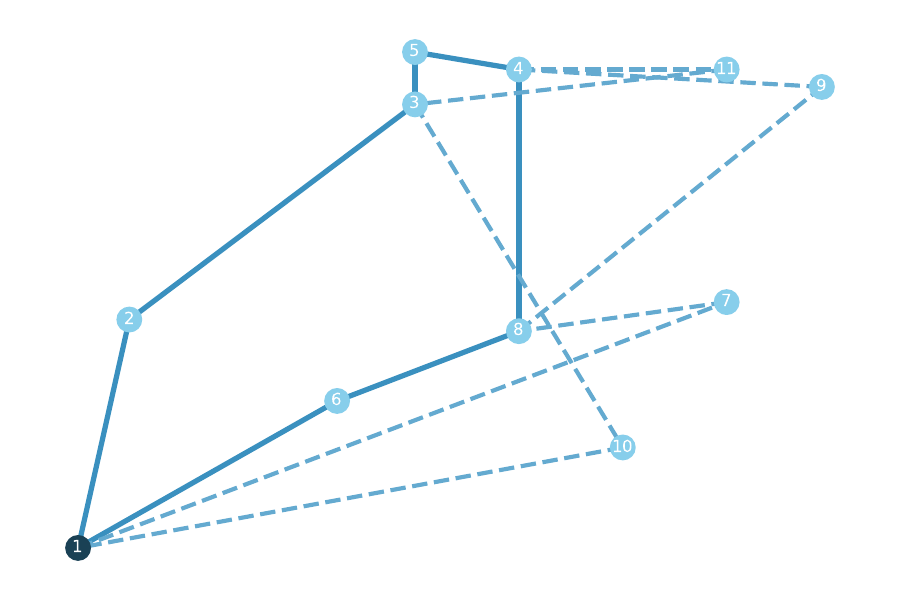}
    Cost = 227.90 (Gap = 0.00\%) \\
    \hline    
    \end{tabular}
    \caption{Routing examples: optimal vs. HM greedy (Instances 6 to 10). Node 1 is the depot, the solid line is the truck route, and the dashed line is the drone route.}%
    \label{fig:example2}
\end{figure}

We evaluate HM on 10 instances of 11-node graphs from \citet{agatz2018optimization}, for which optimal solutions are already known. 
Although HM generates solutions sequentially without backtracking, it can find routes that are close to optimal routes.
The greedy decoding from HM shows a 0.93\% optimality gap, while the sampling strategy reduces the optimality gap to 0.69\%. 
Figures~\ref{fig:example1} and \ref{fig:example2} visualize the optimal solutions reported in the literature and the solution produced by HM with greedy decoding. 
\citet{agatz2018optimization} provide 10 instances of size $N=11$ with exact optimal solutions, available at
the `TSP-D-Instances' repository: \url{https://github.com/pcbouman-eur/TSP-D-Instances}.
We provide the performance of HM greedy and HM sampling on these benchmark instances in Table \ref{tb:optimal}. On average, the optimality gap of HM greedy and HM sampling is 0.93\% and 0.69\%, respectively.

The detailed routes are provided in Figures \ref{fig:example1} and \ref{fig:example2}. 
Note that in Instances 6 and 10, the two routes are not identical, but both are optimal. 
This is because the objective in TSP-D is to minimize the makespan, not the total distance traveled.

\begin{table}
	\caption{The objective values of the optimal solutions and Hybrid Model solutions on TSP-D with 11 nodes.}
	\label{tb:optimal}
	\centering 
	\begin{tabular}{r r r r}
		\toprule
		Instance & Optimal & HM (greedy) (Gap) & HM (1200) (Gap) \\ \midrule
		1        & 221.19  & 223.41 (1.00\%) & 223.41 (1.00\%)        \\
		2        & 205.76  & 205.76 (0.00\%) & 205.76 (0.00\%)         \\
		3        & 192.96  & 193.99 (0.53\%) & 193.99 (0.53\%)         \\
		4        & 241.26  & 241.26 (0.00\%) & 241.26 (0.00\%)         \\
		5        & 248.14  & 249.85 (0.69\%) & 248.82 (0.67\%)         \\
		6        & 217.69  & 217.69 (0.00\%) & 217.69 (0.00\%)         \\
		7        & 237.34  & 240.47 (1.32\%) & 237.34 (0.00\%)         \\
		8        & 214.77  & 226.64 (5.53\%) & 225.36 (4.93\%)         \\
		9        & 256.34  & 256.83 (0.19\%) & 256.83 (0.19\%)         \\
		10       & 227.90  & 227.90 (0.00\%) & 227.90 (0.00\%)         \\ \midrule 
		{mean}     & {226.33}  & {228.38} {(0.93\%)}       &{ 227.84} {(0.69\%)}         \\ \bottomrule
	\end{tabular}
\end{table}

\begin{table*}
    \centering
	\caption{TSP-D results on Random locations dataset. 
	Averages and standard deviation of 100 problem instances. 
	`Cost' refers to the average cost value \eqref{eq:cost}, where the small numbers represent the standard deviation. 
	`Gap' is the mean relative difference to the cost of the best algorithm for each instance as in \eqref{gap}. 
	`Time' is the average solution time of the algorithm for a single instance. 
	}
    \label{tab:TSPD_uniform}
    \begin{adjustbox}{max width=\textwidth}
        \begin{tabular}{l crr crr crr }    
        \toprule
             & \multicolumn{3}{c}{$N=20$} 
             & \multicolumn{3}{c}{$N=50$} 
             & \multicolumn{3}{c}{$N=100$} \\
         \cmidrule(lr){2-4} 
         \cmidrule(lr){5-7} 
         \cmidrule(lr){8-10} 
             Method &
             Cost & Gap & Time & 
             Cost & Gap & Time & 
             Cost & Gap & Time \\
         \midrule
               TSP-ep-all   %
    						&     281.62{\tiny$\pm$18.05} &     0.64\% & (0.17s)
    						&     397.21{\tiny$\pm$20.19} &     1.43\% & (43s)
    						& \bf 535.50{\tiny$\pm$21.83} & \bf 0.78\% & (3992s) \\
               DPS/10       %
    						&     292.23{\tiny$\pm$19.00} &     4.48\% & (0.02s)
    						&     420.51{\tiny$\pm$23.98} &     7.39\% & (0.08s)
    						&     570.74{\tiny$\pm$20.61} &     7.55\% & (0.35s) \\
               DPS/25       %
    						&          - &          - &
    						&     404.78{\tiny$\pm$22.03} &     3.37\% & (1.04s)
    						&     548.23{\tiny$\pm$22.36} &     3.19\% & (2.27s) \\   
		\midrule
               AM (greedy)  %
    						&     294.88{\tiny$\pm$24.44} &     5.36\% & (0.18s)
    						&     439.21{\tiny$\pm$28.09} &    12.16\% & (0.48s)
    						&     642.82{\tiny$\pm$25.12} &    21.08\% & (0.74s) \\
    		    NM (greedy) &     304.70{\tiny$\pm$22.74} &    8.89\% &(0.20s) 
    		                &     445.44{\tiny$\pm$24.51} &   13.78\% &(0.30s) 
    		                & -  & -  & -  \\
    		
               HM (greedy)  %
    						&     285.59{\tiny$\pm$18.10} &     2.07\% & (0.18s)
    						&     408.51{\tiny$\pm$18.78} &     4.33\% & (0.49s)
    						&     564.42{\tiny$\pm$22.03} &     6.27\% & (1.07s) \\
        \midrule    			
               AM (4800)    %
     						&     285.06{\tiny$\pm$18.14} &     1.54\% & (0.30s)
    						&     411.50{\tiny$\pm$20.02} &     4.03\% & (0.78s)
    						&     602.95{\tiny$\pm$19.92} &   13.55\% & (1.92s)\\
               HM (100)     %
    						&     282.65{\tiny$\pm$17.71} &     1.01\% & (0.26s)
    						&     399.11{\tiny$\pm$17.59} &     1.93\% & (0.61s)
    						&     550.01{\tiny$\pm$19.49} &     3.56\% & (1.23s) \\
               HM (1200)    %
    						&     281.90{\tiny$\pm$17.45} &     0.75\% & (0.29s)
    						&     396.80{\tiny$\pm$17.10} &     1.35\% & (0.88s)
    						&     546.28{\tiny$\pm$18.99} &     2.86\% & (1.95s) \\
               HM (2400)    %
    						&     281.75{\tiny$\pm$17.47} &     0.70\% & (0.30s)
    						&     396.16{\tiny$\pm$17.08} &     1.18\% & (1.14s)
    						&     545.32{\tiny$\pm$19.02} &     2.68\% & (2.58s) \\
               HM (4800)    %
    						& \bf 281.53{\tiny$\pm$17.35} & \bf 0.62\% & (0.40s)
    						& \bf 395.65{\tiny$\pm$17.00} & \bf 1.05\% & (1.72s)
    						&     544.55{\tiny$\pm$19.01} &     2.53\% & (4.09s)\\
         \bottomrule
        \end{tabular}
    \end{adjustbox}
\end{table*}

\begin{figure}
    \centering
    \includegraphics[width=0.45\textwidth]{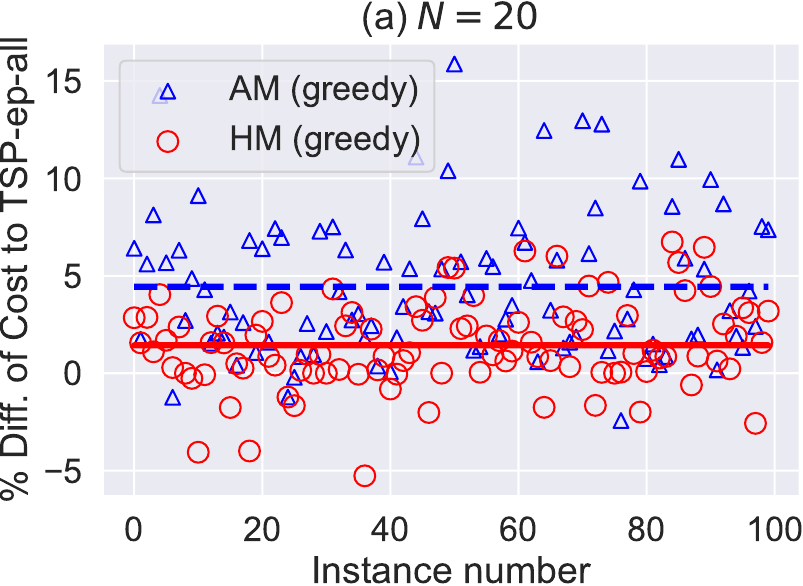}
    \includegraphics[width=0.45\textwidth]{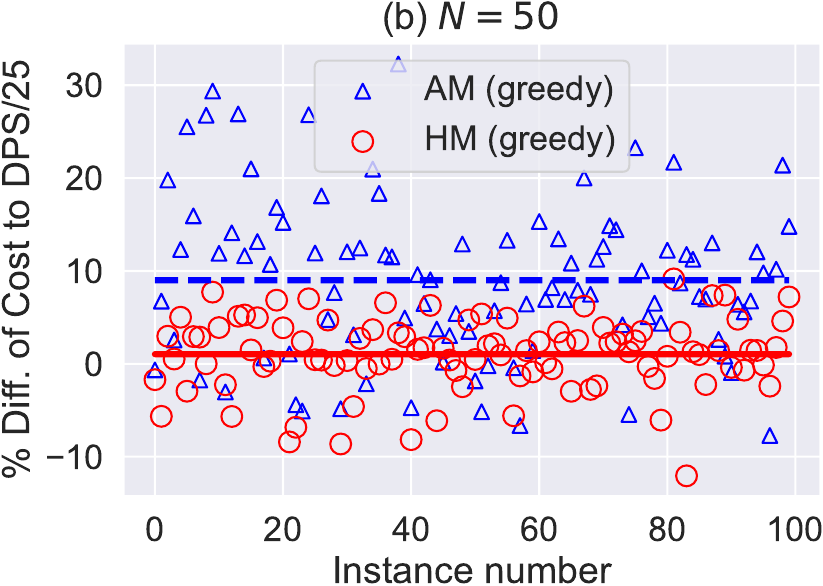}
    \includegraphics[width=0.45\textwidth]{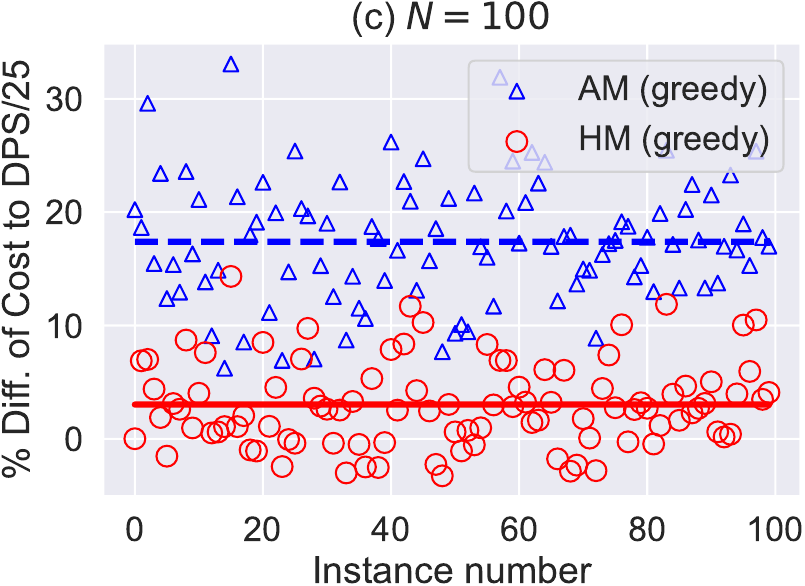}
\caption{AM (greedy) vs. HM (greedy). The percentage difference to TSP-ep-all ($N=20$) or DPS/25 ($N=50, 100$) is reported for each instance. The dashed and solid lines represent the mean values of AM and HM, respectively. }
    \label{fig:AM_HM}
\end{figure}

\begin{figure}
    \centering
	\includegraphics[width=0.6\textwidth]{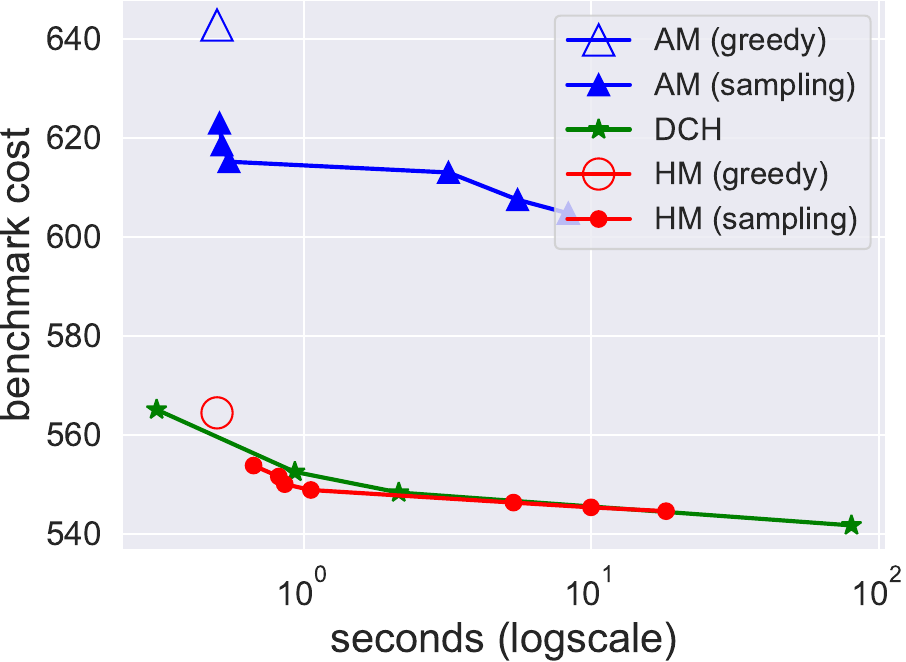}
    \caption{DPS/$g$ vs. AM ($s)$ vs. HM ($s$) for various $g$ and $s$ values. The average cost to the solution time is reported.}
    \label{fig:AM_HM_DPS}
\end{figure}

HM outperforms AM and NM consistently by a large margin, and this trend signifies as the problem size increases.
For example, for $N=100$ in the random location datasets shown in Table \ref{tab:TSPD_uniform}, HM (greedy) has a gap of 6.27\%, while AM (greedy) has a gap of 21.08\% when compared to the best performing method. 
Since NM performed even worse than AM for $N=20$ and $N=50$, we did not test NM for $N=100$.
Figures~\ref{fig:AM_HM}(a)--(c) show the increasing advantage of HM (greedy) over AM (greedy) as $N$ grows.
This demonstrates that recurrent hidden states we introduce in HM are crucial for complex coordination between vehicles that larger problem instances require.

\begin{table*}
    \centering
	\caption{TSP-D results on the Amsterdam dataset. Averages of 100 problem instances.} 
    \label{tab:TSPD_Amsterdam}
    \begin{adjustbox}{max width=\textwidth} 
    \begin{tabular}{l crr crr crr}
    \toprule
         & \multicolumn{3}{c}{$N=10$} 
         & \multicolumn{3}{c}{$N=20$} 
         & \multicolumn{3}{c}{$N=50$} \\
     \cmidrule(lr){2-4} 
     \cmidrule(lr){5-7} 
     \cmidrule(lr){8-10} 
           Method &
         Cost & Gap & Time & Cost & Gap & Time & Cost & Gap & Time \\
     \midrule
           TSP-ep-all   & \bf 2.02{\tiny$\pm$0.23} & \bf  0.69 \% & {(0.01s)}  %
						& \bf 2.35{\tiny$\pm$0.21} & \bf  0.99 \% & {(0.15s)}  %
						& \bf 3.26{\tiny$\pm$0.22} & \bf  0.83 \% & {(42s)} \\ %
		   DPS/10       &     -    &     -       & 
						&     2.49{\tiny$\pm$0.24} &      6.82 \% & {(0.03s)} %
						&     3.55{\tiny$\pm$0.25} &      9.90 \% & {(0.11s)} \\ %
           DPS/25       &     -    &     -       & 
						&     -    &     -       & 
						&     3.37{\tiny$\pm$0.22} &      4.14 \% & {(1.04s)} \\ %

    \midrule
           AM (greedy)  &     2.16{\tiny$\pm$0.29} &      7.38 \% & (0.12s) %
         				&     2.61{\tiny$\pm$0.32} &     11.95 \% & (0.18s) %
						&     4.06{\tiny$\pm$0.43} &     25.24 \% & (0.47s) \\ %
           HM (greedy)  &     2.08{\tiny$\pm$0.26} &      3.46 \% & (0.13s) %
						&     2.49{\tiny$\pm$0.26} &      6.97 \% & (0.19s) %
						&     3.50{\tiny$\pm$0.28} &      8.15 \% & (0.49s) \\ %
    \midrule   
           AM (4800)    &     2.08{\tiny$\pm$0.25} &      3.63 \% & (0.18s) %
						&     2.44{\tiny$\pm$0.23} &      4.62 \% & (0.27s) %
						&     3.72{\tiny$\pm$0.34} &     14.85 \% & (0.89s) \\ %
           HM (100)     &     2.05{\tiny$\pm$0.25} &      2.14 \% & (0.15s) %
						&     2.41{\tiny$\pm$0.24} &      3.27 \% & (0.28s) %
						&     3.36{\tiny$\pm$0.25} &      3.84 \% & (0.65s) \\ %
           HM (1200)    &     2.05{\tiny$\pm$0.25} &      1.75 \% & (0.22s) %
						&     2.39{\tiny$\pm$0.23} &      2.53 \% & (0.33s) %
						&     3.33{\tiny$\pm$0.24} &      2.81 \% & (0.82s) \\ %
		   HM (2400)    &     2.04{\tiny$\pm$0.24} &      1.56 \% & (0.24s) %
						&     2.39{\tiny$\pm$0.23} &      2.36 \% & (0.39s) %
						&     3.32{\tiny$\pm$0.24} &      2.55 \% & (0.99s) \\ %
           HM (4800)    &     2.04{\tiny$\pm$0.24} &      1.58 \% & (0.32s) %
						&     2.38{\tiny$\pm$0.23} &      2.24 \% & (0.52s) %
						&     3.31{\tiny$\pm$0.24} &      2.40 \% & (1.41s) \\ %
     \bottomrule
    \end{tabular}
    \end{adjustbox}
\end{table*}

HM provides superb scalability as well as solid performance, comparable to or exceeding the best heuristic method.
While the performance and speed of HM (greedy) lie between DPS/10 and DPS/25, the batch sampling approach using HM produces the best performance, as shown in Table \ref{tab:TSPD_uniform}, without increasing the solution time as much as TSP-ep-all.
HM's scalability is comparable to DPS/25, but HM outperformed DPS/25 in terms of the average cost. 
Figure~\ref{fig:AM_HM_DPS} also illustrates HM (sampling) is more effective than DPS and AM (sampling).
No direct comparison of solution time between DPS and HM (sampling) is possible, although HM (sampling) shows high efficiency in our experiments.
As DPS can be also run in parallel on multiple CPU threads, the speed of DPS can be further boosted.
However, such parallelization for DPS is not as straightforward as in batch sampling on GPU.

The above observations are consistent in the real Amsterdam dataset in Table~ \ref{tab:TSPD_Amsterdam}: HM provides better solutions than AM.
For $N=50$, AM (greedy) and HM (greedy) had gap of 25.24\% and 8.15 \%, respectively, compared to the TSP-ep-all solutions. 
AM (sampling) and HM (sampling) show the same pattern. 
However, the best solutions are found by TSP-ep-all in the Amsterdam dataset, while HM (sampling) still outperforms DPS.
When customer locations follow non-uniform distributions, the performances of AM and HM seem to deteriorate, which is consistent with the finding of \citet{bogyrbayeva2021reinforcement}.

\begin{table*}
    \centering
	\caption{TSP-D results on the instances of the `TSP-D-Instances' repository (uniform distribution, $\alpha=2$).
The computational time for HGVNS (marked by *) is as reported in \citet{de2020variable}.
	}
    \label{tab:TSPD_hgvns}
    \begin{adjustbox}{max width=\textwidth}
        \begin{tabular}{l crr crr crr }    
        \toprule
             & \multicolumn{3}{c}{$N=20$} 
             & \multicolumn{3}{c}{$N=50$} 
             & \multicolumn{3}{c}{$N=100$} \\
         \cmidrule(lr){2-4} 
         \cmidrule(lr){5-7} 
         \cmidrule(lr){8-10} 
             Method &
             Cost & Gap & Time & 
             Cost & Gap & Time & 
             Cost & Gap & Time \\
         \midrule
               TSP-ep-all   
    						& \bf 276.70 & {\bf 0.85\%} &    {(0.2s)}
    						& \bf 408.16 & {\bf 0.79\%} &    {(40s)}
    						& \bf 548.09 & {\bf 0.63\%} &    {(3635s)} \\
               DPS/25       
    						&          - & -      &    
    						&     421.94 & {4.23\%} &    {(1.1s)}
    						&     559.58 & {2.74\%} &    {(2.0s)} \\   
    	\midrule
               HGVNS         
                            &    293.60 &  &     (0.9s*)
                            &    420.80 &  &     (2.3s*)
                            &    553.43 &  &     (37.8s*) \\
		\midrule
               HM (greedy)  
    						&     277.87 & {1.40\%} &    (0.2s)
    						&     426.93 & {5.53\%} &    (0.4s)
    						&     579.59 & {6.48\%} &    (0.8s) \\
               HM (4800)    
    						&     276.09 & {0.73\%} &    (1.0s)
    						&     408.67 & {1.08\%} &    (3.9s)
    						&     556.82 & {2.20\%} &    (14.4s) \\
         \bottomrule
        \end{tabular}
    \end{adjustbox}
\end{table*}

In Table \ref{tab:TSPD_hgvns}, we also report the performance of HM on larger instances ($N=20,50,100$) from the `TSP-D-Instances' repository, in particular the instances generated by uniform distributions. As before, we assumed $\alpha=2$. 
Note that these instances are generated using the identical method as the random data set in Table \ref{tab:TSPD_uniform}, but there are only 10 instances for each size instead of 100 instances.
We compare the results with the performance of the hybrid general variable neighborhood search (HGVNS) method as reported in \citet{de2020variable} where the algorithm is implemented in C++ and the computational time (marked by *) is measured using the Intel i7 3.6 GHz CPU and 16GB RAM with Ubuntu.
We observe that HM (4800) finds solutions similar to TSP-ep-all for $N=20$ and $N=50$, while it does not perform as good for $N=100$, on average. 
On the other hand, HGVNS finds good solutions for $N=100$, but not as good for $N=20$ and $N=50$, on average.
Since \citet{de2020variable} reported the average cost values only, the mean gap values in Table \ref{tab:TSPD_hgvns} are based on the solutions by other algorithms.

\subsection{Results on TSP-D with Limited Flying Range}

The presented HM can be easily adapted to solve TSP-D with the limited flying range of the drone. 
In this case, the state of MDP also includes the remaining battery of the drone at time $t$. 
Thus, we have $\vec{s}_t := (\Vc_t, c_t^{\tr}, c_t^{\dr}, r_t^{\tr}, r_t^{\dr}, \beta_t^{\dr})$, where 
$\beta_t^{\dr}$ is the battery of the drone at time $t$. We update the battery as follows $\beta_{t+1}^{\dr} = \beta_t^{\dr} - C_t$ if the drone and truck do not travel together and the drone does not wait at its current location. Also, we set the battery to its maximum level when the truck retrieves the drone. 
Next, we pass the state information as an input to the decoder at each step.
In particular, after element-wise projection of the current battery, we use this projection along with the hidden state of the LSTM, the graph embedding, and the traveling time from the current location of the drone to compute the attention vector similar to \eqref{eq:att_vec} as follows:
\begin{equation} \label{eq:att_vec2}
	\vec{a}_{i, \boldsymbol{\cdot}} = \vec{v}_a^\top \tanh \bigg(\mat{W}^a [\bar{\vec{h}};  \vec{\tilde{h}}_{t+1}; \mat{W}^d \vec{\tau}_{i,\boldsymbol{\cdot}}, \mat{W}^{\beta}\vec{\beta_t}] \bigg),
\end{equation}

We also note that the current battery of the drone is passed to the decoder only when selecting the next node to be visited by the drone. 
In order to speed up the training, we also use the current battery level of the drone for masking. For instance, any unserved nodes located further than the current battery level of the drone cannot be selected as the next action. Further, we also make sure to check that there is at least one unserved node where the drone can fly after serving such a node, making sure the drone will have enough battery to be retrieved by the truck.

Table \ref{tab:TSPD_fly_range} presents the performance comparison of TSP-ep-all and HM on benchmark instances from `TSP-D-Instances' repository with uniform distribution and restricted maximum radios for the drone's flying range. We used 60\%, 20\%, and 20\% to set maximum flying ranges for graphs with 10, 20, and 50 nodes, respectively. As shown in Table 5, HM produces comparable results with TSP-ep-all heuristic in solution quality for small-sized instances. %

\begin{table}[]
\centering 
\caption{TSP-D results with limited flying range for drone.} 
    \label{tab:TSPD_fly_range}
    \begin{adjustbox}{max width=\textwidth} 
    \begin{tabular}{l crr crr crr}
    \toprule
         & \multicolumn{3}{c}{$N=10$} 
         & \multicolumn{3}{c}{$N=20$} 
         & \multicolumn{3}{c}{$N=50$} \\
     \cmidrule(lr){2-4} 
     \cmidrule(lr){5-7} 
     \cmidrule(lr){8-10} 
           Method &
         Cost & Gap & Time & Cost & Gap & Time & Cost & Gap & Time \\
     \midrule
TSP-ep-all  & \bf 277.12{\tiny$\pm$28.75}   & \bf 0.00\%   &  {(0.01s)}
            & \bf 387.57{\tiny$\pm$28.19}   & \bf 0.00\%   &  {(0.05s)}      
            & \bf 514.65{\tiny$\pm$25.85}   & \bf 0.00\%   &  {(2.79s)}       \\
\midrule
HM (greedy) & 279.32{\tiny$\pm$27.53}   & 0.85\%   &  (0.07s)      
            & 395.35{\tiny$\pm$33.38}   & 1.94\%   &  (0.09s)     
            & 550.31{\tiny$\pm$32.76}   & 6.93\%   &  (0.21s)      \\
\midrule
HM (100)    & 278.61{\tiny$\pm$28.29}   & 0.56\%   &  (0.21s)  
            & 393.65{\tiny$\pm$32.29}   & 1.51\%   &  (0.40s)  
            & 535.04{\tiny$\pm$29.79}   & 3.96\%   &  (1.05s)   \\
HM (1200)   & 278.18{\tiny$\pm$28.26}   & 0.40\%   &  (0.21s)  
            & 393.00{\tiny$\pm$31.92}   & 1.36\%   &  (0.36s)  
            & 531.39{\tiny$\pm$29.62}   & 3.25\%   &  (1.39s)   \\
HM (2400)   & 278.13{\tiny$\pm$28.29}   & 0.38\%   &  (0.20s)  
            & 393.09{\tiny$\pm$32.05}   & 1.37\%   &  (0.43s)  
            & 530.98{\tiny$\pm$28.75}   & 3.18\%   &  (1.84s)   \\
HM (4800)   & 278.15{\tiny$\pm$28.26}   & 0.40\%   &  (0.27s)  
            & 392.93{\tiny$\pm$32.10}   & 1.33\%   &  (0.57s)  
            & 529.97{\tiny$\pm$29.71}   & 2.97\%   & (3.23s)  \\ \bottomrule
    \end{tabular}
    \end{adjustbox}
\end{table}

\section{Additional Experiments}
\label{sec:additional}

In this section, we provide results from additional experiments to show the performance of the proposed neural network structure and the training algorithm.
First, we demonstrate the performance of the proposed DFPG training algorithm with A2C, when \emph{revisiting nodes} is allowed.
Second, we show that HM is also effective for another class of routing problems, the min-max Capacitated Vehicle Routing Problem (mmCVRP).
Third, we test if HM can perform well on the classical TSP instances.

\subsection{Distributed Training Algorithms for TSP-D when Revisiting is Allowed} \label{sec:revisiting}

In Instance 9 of Figure~\ref{fig:example2}, the exact optimal solution revisits node 9. 
This is an important observation made by \citet{agatz2018optimization}. 
To learn the true optimal solutions, we need to allow revisiting nodes in the simulator.

\begin{figure}
    \centering
	\includegraphics[scale=0.5, keepaspectratio]{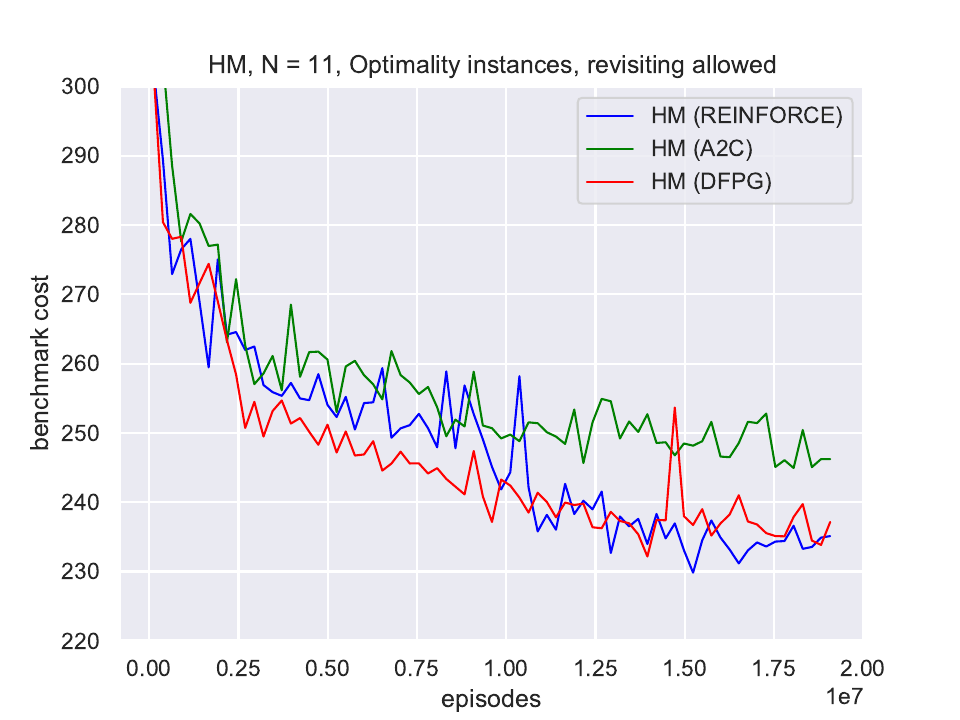}
	\includegraphics[scale=0.5, keepaspectratio]{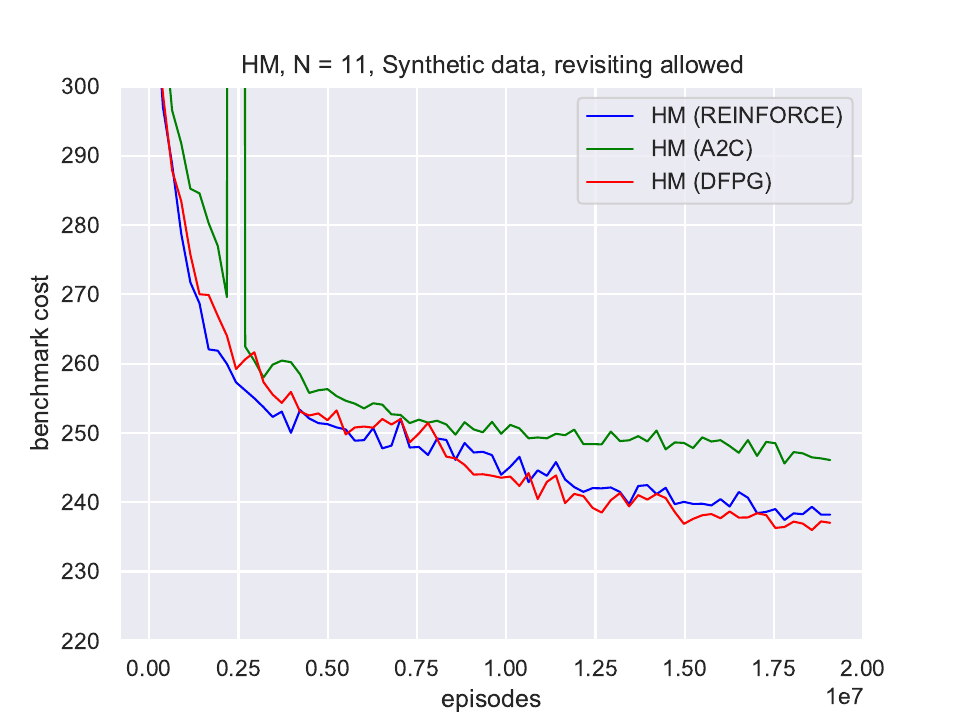}
	\includegraphics[scale=0.5, keepaspectratio]{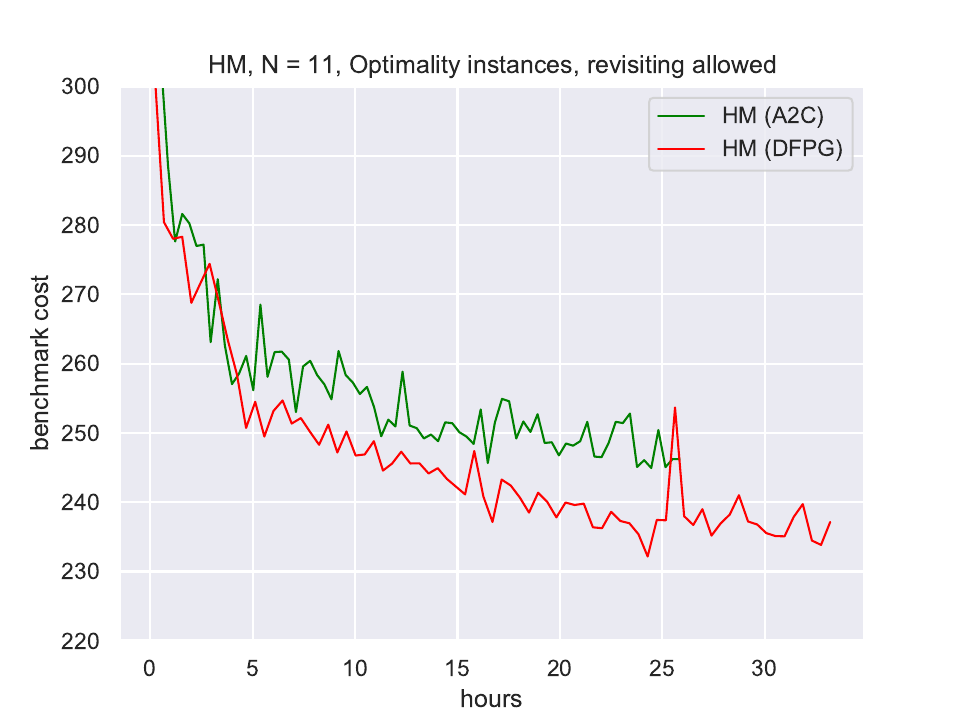}
	\includegraphics[scale=0.5, keepaspectratio]{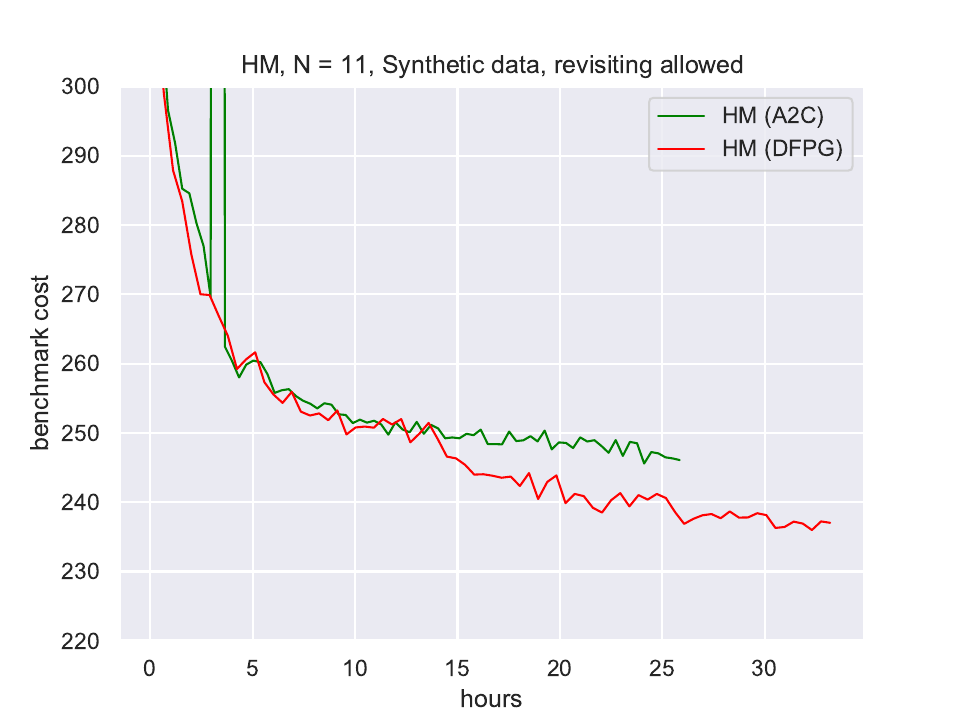}
	\caption{The average benchmark cost curves of HM with REINFORCE, A2C, and DFPG on 11-node graphs from the environment that \emph{allows revisiting}. Optimality instances refer to 11-node instances from \citet{agatz2018optimization}. Synthetic data refers to 100 benchmark instances used in Table~\ref{tab:TSPD_uniform}. The top figures represent the average benchmark costs over the number of episodes. Bottom figures compare the average benchmark costs over wall-clock time.}
	\label{fig:env_revisit_n11}
\end{figure}

Allowing revisiting nodes makes training of both AM and HM significantly slower, as it induces much larger action spaces. AM even did not train smoothly in our experiments. 
When revisiting nodes is allowed, as in Instance 9 in Figure~\ref{fig:example2}, we compared the performances of training algorithms in Figure~\ref{fig:env_revisit_n11}. 
We observe that there are no significant gains from using A2C in the environment with revisiting allowed. 
The proposed DFPG method shows more stable learning curves with the best performance. 
This result demonstrates that DFPG has the potential to be useful in various routing problems with very large action spaces. 

However, in this paper, we did not allow revisiting nodes because our experiments showed not much gain in the performance on average when compared to the results with no revisiting allowed.
For example, the HM solution for Instance 9 without revisiting has a gap of 0.19\% only. Table \ref{tb:env_revisit_opt} shows the testing results of HM greedy and sampling on the Agatz's instances with 11 nodes, which in general is worse than the results with the environment not allowing revisiting nodes reported in Table \ref{tb:optimal}.

\begin{table}[]
\caption{The HM performance when revisiting the nodes is allowed on the instances from Table \ref{tb:optimal}}
\label{tb:env_revisit_opt}
\centering
\begin{tabular}{@{}rrr@{}}
\toprule
Instance & HM (greedy) (Gap) & HM (1200) (Gap) \\ \midrule
1        & 223.41 ( 1.00\%)   & 223.41 ( 1.00\%)     \\
2        & 206.82 ( 0.52\%)   & 206.82 ( 0.52\%)     \\
3        & 214.46 (11.14\%)   & 195.70 ( 1.42\%)     \\
4        & 241.26 ( 0.00\%)   & 241.26 ( 0.00\%)     \\
5        & 257.91 ( 3.94\%)   & 248.82 ( 0.27\%)     \\
6        & 218.37 ( 0.31\%)   & 218.37 ( 0.31\%)     \\
7        & 246.53 ( 3.87\%)   & 237.34 ( 0.00\%)     \\
8        & 223.88 ( 4.24\%)   & 223.88 ( 4.24\%)     \\
9        & 261.12 ( 1.86\%)   & 257.02 ( 0.27\%)     \\
10       & 227.90 ( 0.00\%)   & 227.90 ( 0.00\%)     \\ \midrule
mean     & 232.17 ( 2.69\%)   & 228.05 ( 0.80\%)     \\ \bottomrule
\end{tabular}
\end{table}

\subsection{Additional Results on mmCVRP}
While our main focus is TSP-D, we also test the performances of HM and AM in mmCVRP to demonstrate that HM can generalize to other multi-vehicle routing problems. 
TSP-D concerns the routing of heterogeneous vehicles (the drone is faster than the truck), whereas mmCVRP is about the routing of homogeneous vehicles (all vehicles are at the same speed and of the same capacity). 
In both cases, we aim to minimize the makespan to finish all jobs.
The details of the experiments follow.

\paragraph{Problem statement} \label{sec:mmCVRP}
We have a set of nodes representing customer locations and a depot. 
Each customer has a known demand, which a fleet of vehicles must satisfy and return back to a depot in a minimal time. 
Each customer can be served by exactly one vehicle. 
The routing decisions of vehicles are restricted by their capacity to carry the load.

We compare HM against OR Tools \citep{ortools} to solve mmCVRP. 
We implemented the mmCVRP simulator using Julia 1.5 and interfaced it with Python using pyjulia available at \url{https://github.com/JuliaPy/pyjulia}.
For both $N=20$ and $N=30$, we used three vehicles. 
We use the same set of hyperparameters for both TSP-D and mmCVRP.

\paragraph{Data generation}
We generate the random data for mmCVRP by sampling the x and y coordinates of the customer locations using Uniform(0, 1).
We randomly assign customer demand at each node with integer values ranging between 1 and 9.
We set the capacity of each vehicle as follows:
\begin{equation}
    Q = \bigg\lceil 1.2/m\sum_{n \in \Nc}D_n \bigg\rceil, 
\end{equation}
where $D_i$ is demand at node $i$, and $m$ is the number of vehicles.
We computed the traveling time between nodes as follows:
\begin{align}
    \tau_{i, j} = \Big\lfloor 2000d_{ij} \Big\rfloor \quad  \forall i \in \Nc, j \in \Nc,
\end{align}
where $d_{i,j}$ is the Euclidean distance between nodes $i$ and $j$. 

\paragraph{HM implementation for mmCVRP} 
We embed a graph exactly as in TSP-D. However, in the decoder, we also pass information about the remaining capacity of a vehicle. In particular, for each node $i$, we compute expected capacity if it is visited by a decision taker at time $t$: $Q'_t= Q_t-D^t_i$ for all $i \in N$, where $Q_t$ is a capacity of a decision taker at time $t$ and $D^t_i$ is demand at node $i$ at time $t$. Then we use a vector of expected remaining capacities $\vec{Q}'_t$ to compute attention in the decoder:

\begin{equation}
	\vec{a}_{i, \boldsymbol{\cdot}} = \vec{v}_a^\top \tanh \bigg(\mat{W}^a [\bar{\vec{h}};  \vec{\tilde{h}}_{t+1}; \mat{W}^d \vec{Q'}_t] \bigg),
\end{equation}

\paragraph{AM implementation for mmCVRP} 
To compute the context node for mmCVRP, we use the embeddings of the current node of a decision taker and embeddings of the nodes selected by other nodes. In particular, in the presence of $m$ vehicles in mmCVRP, we will have:

\begin{equation}
   \vec{h}_c = [ \vec{\bar{h}}, \vec{h}^{L}_{i},\vec{h}^{L}_{j_1} ... \vec{h}^{L}_{j_{m-1}}, \mat{W}^d \vec{Q}_t] 
\end{equation}
where $\vec{h}^{L}_{i}$ is embedding of a node $i$ selected by a decision taker and $\vec{h}^{L}_{j_{m-1}}$ is a embedding of a node $j$ selected by a vehicle $m-1$. In computing the context node, we also use the remaining capacity of a decision taker $\vec{Q}_t$ at time $t$.

\begin{table*}[t]
    \centering
	\caption{mmCVRP results with three vehicles. Averages of 128 problem instances.} 
    \label{tb:mmCVRP}
    \begin{adjustbox}{max width=\textwidth}
    \begin{tabular}{l crr crr}
    \toprule
         & \multicolumn{3}{c}{$N=20$} 
         & \multicolumn{3}{c}{$N=30$} \\
     \cmidrule(lr){2-4} 
     \cmidrule(lr){5-7} 
           Method &
         Cost & Gap & Time & 
         Cost & Gap & Time \\
     \midrule
           OR-Tools     & \bf 4,086.99{\tiny$\pm$471.78} & \bf 0.00 \% & (1.0s)
						&     4,871.06{\tiny$\pm$453.74} &     7.26 \% & (1.0s) \\
           AM (greedy)  &     4,396.11{\tiny$\pm$535.81} &     7.56 \% & (0.1s) %
						&     4,939.37{\tiny$\pm$473.15} &     8.76 \% & (0.2s) \\ %
           HM (greedy)  &     4,298.83{\tiny$\pm$523.19} &     5.18 \% & (0.1s) %
						&     4,779.27{\tiny$\pm$488.44} &     5.24 \% & (0.2s) \\ %
		   AM (1200)    &     4,142.74{\tiny$\pm$478.88} &     1.36 \% & (14.3s) %
						&     4,574.05{\tiny$\pm$420.01} &     0.01 \% & (27.6s) \\ %
           HM (1200)    &     4,129.41{\tiny$\pm$496.54} &     1.04 \% & (11.4s) %
						& \bf 4,541.45{\tiny$\pm$429.89} & \bf 0.00 \% & (25.2s) \\ %
     \bottomrule
    \end{tabular}
    \end{adjustbox}
\end{table*}

\paragraph{Results}
The results in Table \ref{tb:mmCVRP} demonstrate that HM provides competitive results with OR-Tools to solve the mmCVRP, while AM was not as competitive. 
This indicates that the proposed HM can be effective in learning solutions to other coordinated routing problems, especially when the objective is to minimize the makespan.

\subsection{Performances on TSP Instances} \label{sec:tsp}

We compare the performance of AM, HM, and the RL model by \citet{nazari2018reinforcement} for TSP. 
We call the \citeauthor{nazari2018reinforcement} model NM.
We generate random location instances with the same method as \citet{kool2018attention} and \citet{nazari2018reinforcement}. 
While we implemented AM and HM for TSP by ourselves, we used the result reported by \citet{nazari2018reinforcement} for NM. 
Random location instances are generated using the same rule.
\textbf{HM for TSP} is similar to TSP-D, except we pass a distance vector, $\vec{d}_{i,\boldsymbol{\cdot}}$, from the current location of an agent $i$ to other nodes as follows in computing attention in decoder:
\begin{equation*}
	\vec{a}_{i, \boldsymbol{\cdot}} = \vec{v}_a^\top \tanh \bigg(\mat{W}^a [\bar{\vec{h}};  \vec{\tilde{h}}_{t+1}; \mat{W}^d \vec{d}_{i,\boldsymbol{\cdot}}] \bigg).
\end{equation*}

The result in Table~\ref{tab:TSP} shows that HM is not as competitive as AM in the classical TSP instances, while HM performs as similar as NM.

\begin{table*}[h!]
    \centering
	\caption{TSP results on random locations. An average of 10,000 problem instances is reported. Cost refers to the cost function \eqref{eq:cost}. The gap is the relative difference to the cost of the best algorithm for the setting.} 
    \label{tab:TSP}
    \begin{adjustbox}{max width=\textwidth}
    \begin{tabular}{l rr rr rr}
    \toprule
         & \multicolumn{2}{c}{$N=20$} 
         & \multicolumn{2}{c}{$N=50$} 
         & \multicolumn{2}{c}{$N=100$} \\
     \cmidrule(lr){2-3} 
     \cmidrule(lr){4-5} 
     \cmidrule(lr){6-7} 
           Method &
         Cost & Gap & Cost & Gap & Cost & Gap \\
     \midrule
           Concorde     & \bf 3.83 & \bf 0.00 \%
						& \bf 5.69 & \bf 0.00 \%
						& \bf 7.76 & \bf 0.00 \% \\
           AM (greedy)  &     3.84 &     0.29 \% 
         				&     5.79 &     1.67 \% 
						&     8.10 &     4.38 \% \\
           AM (sampling)&     3.83 &     0.07 \% 
         				&     5.72 &     0.48 \% 
						&     7.95 &     2.34 \% \\
           HM (greedy)  &     3.88 &     1.38 \% 
						&     6.00 &     5.39 \% 
						&     8.54 &     9.93 \% \\
           HM (sampling)&     3.84 &     0.39 \% 
						&     5.86 &     2.92 \% 
						&     8.11 &     4.45 \% \\
    \midrule
           NM (greedy)
                        &     3.97 &      
						&     6.08 &     
						&     8.44 &      \\
     \bottomrule
    \end{tabular}
    \end{adjustbox}
\end{table*}

\section{Conclusions} \label{sec:concl}
This study proposed a new end-to-end learning model that can learn the coordinated routing of multiple vehicles. 
Our method applies to both heterogeneous and homogeneous fleet cases, as shown in TSP-D and mmCVRP, respectively.
{Compared to other learning methods to solve routing problems, the proposed model efficiently deals with the coordination of multiple vehicles and achieves comparable performances as strong optimization heuristics methods.}
This could be because the LSTM-based decoder of the proposed model stores the decisions of other vehicles to make better-informed decisions compared to stateless attention-based models. 
Moreover, in contrast to optimization heuristics that solve TSP-D and mmCVRP tailored to the specifics of each problem, the proposed model solves both problems with the same hyperparameters and input structures.
The proposed model, however, does not perform as strongly in TSP.

Searching for a universal architecture that efficiently solves both TSP and TSP-D will be an important direction for future research.
Such a universal architecture will solve a broad class of routing problems arising in various logistics services, including heterogeneous fleets and multiple echelons of service vehicles. 
Identifying architectures and training algorithms robust to variations in distributions of customer locations will be another important direction.
It is unclear how such variations would impact the performance of end-to-end learning methods.
Understanding its impact will be critical for real-world applications.

\section*{Acknowledgments}
This work was supported by the Institute of Information \& Communications Technology Planning \& Evaluation (IITP) grant (2020-0-01336, Artificial Intelligence Graduate School Program, UNIST), 2021 Research Fund (1.210107.01) of UNIST, the National Research Foundation of Korea (2021R1A4A3033149 and 2021H1D3A2A01039401), and the National Science Foundation of the U.S. (CMMI-2032458), the internal grant of Suleyman Demirel University 2021-2022.

\bibliographystyle{ormsv080-ck}
\bibliography{reference.bib}

\begin{thebibliography}{50}
\expandafter\ifx\csname natexlab\endcsname\relax\def\natexlab#1{#1}\fi
\expandafter\ifx\csname url\endcsname\relax
  \def\url#1{{\tt #1}}\fi
\expandafter\ifx\csname urlprefix\endcsname\relax\def\urlprefix{URL }\fi
\expandafter\ifx\csname urlstyle\endcsname\relax
  \expandafter\ifx\csname doi\endcsname\relax
  \def\doi#1{doi:\discretionary{}{}{}#1}\fi \else
  \expandafter\ifx\csname doi\endcsname\relax
  \def\doi{doi:\discretionary{}{}{}\begingroup \urlstyle{rm}\Url}\fi \fi

\bibitem[{Agatz et~al.(2018)Agatz, Bouman, and Schmidt}]{agatz2018optimization}
Agatz, N., P. Bouman, M. Schmidt. 2018.
\newblock Optimization approaches for the traveling salesman problem with
  drone.
\newblock {\it Transportation Science\/} {\bf 52}(4) 965--981.

\bibitem[{Amazon Prime Air(2022)}]{amazon_news}
Amazon Prime Air. 2022.
\newblock Amazon {P}rime {A}ir.
\newblock
  \url{https://www.aboutamazon.com/news/transportation/amazon-prime-air-prepares-for-drone-deliveries}.
\newblock Accessed: August 16, 2022.

\bibitem[{Applegate et~al.(2001)Applegate, Bixby, Chv{\'a}tal, and
  Cook}]{applegate2001tsp}
Applegate, D., R. Bixby, V. Chv{\'a}tal, W. Cook. 2001.
\newblock {TSP} cuts which do not conform to the template paradigm.
\newblock {\it Computational Combinatorial Optimization\/}. Springer, 261--303.

\bibitem[{Bello et~al.(2016)Bello, Pham, Le, Norouzi, and
  Bengio}]{bello2016neural}
Bello, I., H. Pham, Q.~V. Le, M. Norouzi, S. Bengio. 2016.
\newblock Neural combinatorial optimization with reinforcement learning.
\newblock {\it arXiv preprint arXiv:1611.09940\/} .

\bibitem[{Bezanson et~al.(2017)Bezanson, Edelman, Karpinski, and
  Shah}]{bezanson2017julia}
Bezanson, J., A. Edelman, S. Karpinski, V.~B. Shah. 2017.
\newblock Julia: A fresh approach to numerical computing.
\newblock {\it SIAM Review\/} {\bf 59}(1) 65--98.

\bibitem[{Boccia et~al.(2021)Boccia, Masone, Sforza, and
  Sterle}]{boccia2021column}
Boccia, M., A. Masone, A. Sforza, C. Sterle. 2021.
\newblock A column-and-row generation approach for the flying sidekick
  travelling salesman problem.
\newblock {\it Transportation Research Part C: Emerging Technologies\/} {\bf
  124} 102913.

\bibitem[{Bogyrbayeva et~al.(2021)Bogyrbayeva, Jang, Shah, Jang, and
  Kwon}]{bogyrbayeva2021reinforcement}
Bogyrbayeva, A., S. Jang, A. Shah, Y.~J. Jang, C. Kwon. 2021.
\newblock A reinforcement learning approach for rebalancing electric vehicle
  sharing systems.
\newblock {\it IEEE Transactions on Intelligent Transportation Systems\/} {\bf
  Accepted} 1--11.

\bibitem[{Bogyrbayeva et~al.(2022)Bogyrbayeva, Meraliyev, Mustakhov, and
  Dauletbayev}]{bogyrbayeva2022learning}
Bogyrbayeva, A., M. Meraliyev, T. Mustakhov, B. Dauletbayev. 2022.
\newblock Learning to solve vehicle routing problems: A survey.
\newblock {\it arXiv preprint arXiv:2205.02453\/} .

\bibitem[{Bouman et~al.(2018)Bouman, Agatz, and Schmidt}]{bouman2018dynamic}
Bouman, P., N. Agatz, M. Schmidt. 2018.
\newblock Dynamic programming approaches for the traveling salesman problem
  with drone.
\newblock {\it Networks\/} {\bf 72}(4) 528--542.

\bibitem[{Business Insider(2017)}]{UPS}
Business Insider. 2017.
\newblock Ups tests drone delivery system.
\newblock
  \url{https://www.businessinsider.com/ups-tests-drone-delivery-system-2017-2}.
\newblock Accessed: August 18, 2022.

\bibitem[{Carlsson and Song(2018)}]{carlsson2018coordinated}
Carlsson, J.~G., S. Song. 2018.
\newblock Coordinated logistics with a truck and a drone.
\newblock {\it Management Science\/} {\bf 64}(9) 4052--4069.

\bibitem[{Chung et~al.(2020)Chung, Sah, and Lee}]{chung2020optimization}
Chung, S.~H., B. Sah, J. Lee. 2020.
\newblock Optimization for drone and drone-truck combined operations: A review
  of the state of the art and future directions.
\newblock {\it Computers \& Operations Research\/}  105004.

\bibitem[{de~Freitas and Penna(2020)}]{de2020variable}
de~Freitas, J.~C., P.~H.~V. Penna. 2020.
\newblock A variable neighborhood search for flying sidekick traveling salesman
  problem.
\newblock {\it International Transactions in Operational Research\/} {\bf
  27}(1) 267--290.

\bibitem[{Dell'Amico et~al.(2019)Dell'Amico, Montemanni, and Novellani}]{arxiv}
Dell'Amico, M., R. Montemanni, S. Novellani. 2019.
\newblock Models and algorithms for the flying sidekick traveling salesman
  problem.
\newblock \doi{10.48550/ARXIV.1910.02559}.
\newblock \urlprefix\url{https://arxiv.org/abs/1910.02559}.

\bibitem[{Dell'Amico et~al.(2021)Dell'Amico, Montemanni, and
  Novellani}]{dell2021random}
Dell'Amico, M., R. Montemanni, S. Novellani. 2021.
\newblock A random restart local search matheuristic for the flying sidekick
  traveling salesman problem.
\newblock {\it 2021 The 8th International Conference on Industrial Engineering
  and Applications (Europe)\/}. 205--209.

\bibitem[{Dell’Amico et~al.(2021{\natexlab{a}})Dell’Amico, Montemanni, and
  Novellani}]{dell2021algorithms}
Dell’Amico, M., R. Montemanni, S. Novellani. 2021{\natexlab{a}}.
\newblock Algorithms based on branch and bound for the flying sidekick
  traveling salesman problem.
\newblock {\it Omega\/} {\bf 104} 102493.

\bibitem[{Dell’Amico et~al.(2021{\natexlab{b}})Dell’Amico, Montemanni, and
  Novellani}]{dell2021drone}
Dell’Amico, M., R. Montemanni, S. Novellani. 2021{\natexlab{b}}.
\newblock Drone-assisted deliveries: New formulations for the flying sidekick
  traveling salesman problem.
\newblock {\it Optimization Letters\/} {\bf 15}(5) 1617--1648.

\bibitem[{Devlin et~al.(2018)Devlin, Chang, Lee, and
  Toutanova}]{devlin2018bert}
Devlin, J., M.-W. Chang, K. Lee, K. Toutanova. 2018.
\newblock {BERT}: Pre-training of deep bidirectional transformers for language
  understanding.
\newblock {\it arXiv preprint arXiv:1810.04805\/} .

\bibitem[{El-Adle et~al.(2021)El-Adle, Ghoniem, and Haouari}]{el2021parcel}
El-Adle, A.~M., A. Ghoniem, M. Haouari. 2021.
\newblock Parcel delivery by vehicle and drone.
\newblock {\it Journal of the Operational Research Society\/} {\bf 72}(2)
  398--416.

\bibitem[{Gonzalez-R et~al.(2020)Gonzalez-R, Canca, Andrade-Pineda, Calle, and
  Leon-Blanco}]{gonzalez2020truck}
Gonzalez-R, P.~L., D. Canca, J.~L. Andrade-Pineda, M. Calle, J.~M. Leon-Blanco.
  2020.
\newblock Truck-drone team logistics: A heuristic approach to multi-drop route
  planning.
\newblock {\it Transportation Research Part C: Emerging Technologies\/} {\bf
  114} 657--680.

\bibitem[{Ha et~al.(2018)Ha, Deville, Pham, and H{\`a}}]{ha2018min}
Ha, Q.~M., Y. Deville, Q.~D. Pham, M.~H. H{\`a}. 2018.
\newblock On the min-cost traveling salesman problem with drone.
\newblock {\it Transportation Research Part C: Emerging Technologies\/} {\bf
  86} 597--621.

\bibitem[{Haider et~al.(2019)Haider, Charkhgard, Kim, and
  Kwon}]{haider2019optimizing}
Haider, Z., H. Charkhgard, S.~W. Kim, C. Kwon. 2019.
\newblock Optimizing the relocation operations of free-floating electric
  vehicle sharing systems.
\newblock {\it Available at SSRN\/} .

\bibitem[{Hottung and Tierney(2020)}]{hottung2020neural}
Hottung, A., K. Tierney. 2020.
\newblock Neural large neighborhood search for the capacitated vehicle routing
  problem.
\newblock {\it ECAI 2020\/}. IOS Press, 443--450.

\bibitem[{Joerss et~al.(2016)Joerss, Schroder, Neuhaus, Klink, and
  Mann}]{joerss2016parcel}
Joerss, M., J. Schroder, F. Neuhaus, C. Klink, F. Mann. 2016.
\newblock Parcel delivery: The future of last mile.
\newblock Tech. rep., McKinsey \& Company, Travel, Transport and Logistics.
\newblock
  \urlprefix\url{https://www.mckinsey.com/~/media/mckinsey/industries/travel%20logistics%20and%20infrastructure/our%20insights/how%20customer%20demands%20are%20reshaping%20last%20mile%20delivery/parcel_delivery_the_future_of_last_mile.pdf}.

\bibitem[{Khalil et~al.(2017)Khalil, Dai, Zhang, Dilkina, and
  Song}]{dai2017learning}
Khalil, E., H. Dai, Y. Zhang, B. Dilkina, L. Song. 2017.
\newblock Learning combinatorial optimization algorithms over graphs.
\newblock {\it Advances in Neural Information Processing Systems\/} {\bf 30}
  6348--6358.

\bibitem[{Kim et~al.(2021)Kim, Park, and Kim}]{kim2021learning}
Kim, M., J. Park, J. Kim. 2021.
\newblock Learning collaborative policies to solve {NP}-hard routing problems.
\newblock {\it Thirty-Fifth Conference on Neural Information Processing
  Systems\/}.

\bibitem[{Kool et~al.(2021)Kool, van Hoof, Gromicho, and
  Welling}]{kool2021deep}
Kool, W., H. van Hoof, J. Gromicho, M. Welling. 2021.
\newblock Deep policy dynamic programming for vehicle routing problems.
\newblock {\it arXiv preprint arXiv:2102.11756\/} .

\bibitem[{Kool et~al.(2018)Kool, van Hoof, and Welling}]{kool2018attention}
Kool, W., H. van Hoof, M. Welling. 2018.
\newblock Attention, learn to solve routing problems!
\newblock {\it International Conference on Learning Representations\/}.

\bibitem[{Kwon et~al.(2020)Kwon, Choo, Kim, Yoon, Gwon, and Min}]{kwon2020pomo}
Kwon, Y.-D., J. Choo, B. Kim, I. Yoon, Y. Gwon, S. Min. 2020.
\newblock Pomo: Policy optimization with multiple optima for reinforcement
  learning.
\newblock H. Larochelle, M. Ranzato, R. Hadsell, M.~F. Balcan, H. Lin, eds.,
  {\it Advances in Neural Information Processing Systems\/}, vol.~33. Curran
  Associates, Inc., 21188--21198.

\bibitem[{Liu et~al.(2022)Liu, Li, and Khojandi}]{liu2022flying}
Liu, Z., X. Li, A. Khojandi. 2022.
\newblock The flying sidekick traveling salesman problem with stochastic travel
  time: A reinforcement learning approach.
\newblock {\it Transportation Research Part E: Logistics and Transportation
  Review\/} {\bf 164} 102816.

\bibitem[{Lu et~al.(2020)Lu, Zhang, and Yang}]{lu2019learning}
Lu, H., X. Zhang, S. Yang. 2020.
\newblock A learning-based iterative method for solving vehicle routing
  problems.
\newblock {\it International Conference on Learning Representations\/}.
\newblock \urlprefix\url{https://openreview.net/forum?id=BJe1334YDH}.

\bibitem[{Macrina et~al.(2020)Macrina, Pugliese, Guerriero, and
  Laporte}]{macrina2020drone}
Macrina, G., L.~D.~P. Pugliese, F. Guerriero, G. Laporte. 2020.
\newblock Drone-aided routing: A literature review.
\newblock {\it Transportation Research Part C: Emerging Technologies\/} {\bf
  120} 102762.

\bibitem[{Mazyavkina et~al.(2021)Mazyavkina, Sviridov, Ivanov, and
  Burnaev}]{mazyavkina2021reinforcement}
Mazyavkina, N., S. Sviridov, S. Ivanov, E. Burnaev. 2021.
\newblock Reinforcement learning for combinatorial optimization: A survey.
\newblock {\it Computers \& Operations Research\/}  105400.

\bibitem[{Mnih et~al.(2013)Mnih, Kavukcuoglu, Silver, Graves, Antonoglou,
  Wierstra, and Riedmiller}]{mnih2013playing}
Mnih, V., K. Kavukcuoglu, D. Silver, A. Graves, I. Antonoglou, D. Wierstra, M.
  Riedmiller. 2013.
\newblock Playing atari with deep reinforcement learning.
\newblock {\it arXiv preprint arXiv:1312.5602\/} .

\bibitem[{Murray and Chu(2015)}]{murray2015flying}
Murray, C.~C., A.~G. Chu. 2015.
\newblock The flying sidekick traveling salesman problem: Optimization of
  drone-assisted parcel delivery.
\newblock {\it Transportation Research Part C: Emerging Technologies\/} {\bf
  54} 86--109.

\bibitem[{Nazari et~al.(2018)Nazari, Oroojlooy, Snyder, and
  Tak{\'a}c}]{nazari2018reinforcement}
Nazari, M., A. Oroojlooy, L. Snyder, M. Tak{\'a}c. 2018.
\newblock Reinforcement learning for solving the vehicle routing problem.
\newblock {\it Advances in Neural Information Processing Systems\/}.
  9839--9849.

\bibitem[{Perron and Furnon(2019)}]{ortools}
Perron, L., V. Furnon. 2019.
\newblock {OR-Tools}.
\newblock \urlprefix\url{https://developers.google.com/optimization/}.

\bibitem[{Poikonen et~al.(2019)Poikonen, Golden, and
  Wasil}]{poikonen2019branch}
Poikonen, S., B. Golden, E.~A. Wasil. 2019.
\newblock A branch-and-bound approach to the traveling salesman problem with a
  drone.
\newblock {\it INFORMS Journal on Computing\/} {\bf 31}(2) 335--346.

\bibitem[{Roberti and Ruthmair(2021)}]{roberti2021exact}
Roberti, R., M. Ruthmair. 2021.
\newblock Exact methods for the traveling salesman problem with drone.
\newblock {\it Transportation Science\/} {\bf 55}(2) 315--335.

\bibitem[{Schermer et~al.(2020)Schermer, Moeini, and Wendt}]{schermer2020b}
Schermer, D., M. Moeini, O. Wendt. 2020.
\newblock A b ranch-and-cut approach and alternative formulations for the
  traveling salesman problem with drone.
\newblock {\it Networks\/} {\bf 76}(2) 164--186.

\bibitem[{Scott(2015)}]{scott2015multivariate}
Scott, D.~W. 2015.
\newblock {\it Multivariate Density Estimation: Theory, Practice, and
  Visualization\/}.
\newblock John Wiley \& Sons.

\bibitem[{Sykora et~al.(2020)Sykora, Ren, and Urtasun}]{sykora2020multi}
Sykora, Q., M. Ren, R. Urtasun. 2020.
\newblock Multi-agent routing value iteration network.
\newblock {\it International Conference on Machine Learning\/}. PMLR,
  9300--9310.

\bibitem[{V{\'a}squez et~al.(2021)V{\'a}squez, Angulo, and
  Klapp}]{vasquez2021exact}
V{\'a}squez, S.~A., G. Angulo, M.~A. Klapp. 2021.
\newblock An exact solution method for the tsp with drone based on
  decomposition.
\newblock {\it Computers \& Operations Research\/} {\bf 127} 105127.

\bibitem[{Vaswani et~al.(2017)Vaswani, Shazeer, Parmar, Uszkoreit, Jones,
  Gomez, Kaiser, and Polosukhin}]{vaswani2017attention}
Vaswani, A., N. Shazeer, N. Parmar, J. Uszkoreit, L. Jones, A.~N. Gomez, L.
  Kaiser, I. Polosukhin. 2017.
\newblock Attention is all you need.
\newblock {\it arXiv preprint arXiv:1706.03762\/} .

\bibitem[{Vinyals et~al.(2015)Vinyals, Fortunato, and
  Jaitly}]{vinyals2015pointer}
Vinyals, O., M. Fortunato, N. Jaitly. 2015.
\newblock Pointer networks.
\newblock {\it arXiv preprint arXiv:1506.03134\/} .

\bibitem[{Virtanen et~al.(2020)Virtanen, Gommers, Oliphant, Haberland, Reddy,
  Cournapeau, Burovski, Peterson, Weckesser, Bright, {van der Walt}, Brett,
  Wilson, Millman, Mayorov, Nelson, Jones, Kern, Larson, Carey, Polat, Feng,
  Moore, {VanderPlas}, Laxalde, Perktold, Cimrman, Henriksen, Quintero, Harris,
  Archibald, Ribeiro, Pedregosa, {van Mulbregt}, and {SciPy 1.0
  Contributors}}]{2020SciPy-NMeth}
Virtanen, P., R. Gommers, T.~E. Oliphant, M. Haberland, T. Reddy, D.
  Cournapeau, E. Burovski, P. Peterson, W. Weckesser, J. Bright, S.~J. {van der
  Walt}, M. Brett, J. Wilson, K.~J. Millman, N. Mayorov, A.~R.~J. Nelson, E.
  Jones, R. Kern, E. Larson, C.~J. Carey, {\.I}. Polat, Y. Feng, E.~W. Moore,
  J. {VanderPlas}, D. Laxalde, J. Perktold, R. Cimrman, I. Henriksen, E.~A.
  Quintero, C.~R. Harris, A.~M. Archibald, A.~H. Ribeiro, F. Pedregosa, P. {van
  Mulbregt}, {SciPy 1.0 Contributors}. 2020.
\newblock {{SciPy} 1.0: Fundamental Algorithms for Scientific Computing in
  Python}.
\newblock {\it Nature Methods\/} {\bf 17} 261--272.
\newblock \doi{10.1038/s41592-019-0686-2}.

\bibitem[{Williams(1992)}]{williams1992simple}
Williams, R.~J. 1992.
\newblock Simple statistical gradient-following algorithms for connectionist
  reinforcement learning.
\newblock {\it Machine Learning\/} {\bf 8}(3-4) 229--256.

\bibitem[{Wing Aviation(2021)}]{wing_news}
Wing Aviation. 2021.
\newblock Wing delivers library books to students in {V}irginia.
\newblock
  \url{https://blog.wing.com/2020/06/wing-delivers-library-books-to-students.html}.
\newblock Accessed: March 10, 2021.

\bibitem[{Yurek and Ozmutlu(2018)}]{yurek2018decomposition}
Yurek, E.~E., H.~C. Ozmutlu. 2018.
\newblock A decomposition-based iterative optimization algorithm for traveling
  salesman problem with drone.
\newblock {\it Transportation Research Part C: Emerging Technologies\/} {\bf
  91} 249--262.

\bibitem[{Zhang et~al.(2020)Zhang, He, Zhang, Lin, and Li}]{zhang2020multi}
Zhang, K., F. He, Z. Zhang, X. Lin, M. Li. 2020.
\newblock Multi-vehicle routing problems with soft time windows: A multi-agent
  reinforcement learning approach.
\newblock {\it Transportation Research Part C: Emerging Technologies\/} {\bf
  121} 102861.

\end{thebibliography}

\end{document}